\input{style/arxiv-general.cfg}
\documentclass[MSNbibl,citesort,number,seceqn,dvips]{arxbj}
\makeatletter
   \@ifpackageloaded{graphicx}{}{\usepackage{graphicx}}
\makeatother
\usepackage{multirow,upgreek}
\usepackage{graphicx}
\usepackage{mathrsfs}

\aid{0}
\volume{21}
\issue{2}
\pubyear{2015}
\firstpage{1089}
\lastpage{1133}
\doi{10.3150/14-BEJ599} 

\makeatletter

\newcommand{\rright}{\right}
\newcommand{\lleft}{\left}
\newcommand{\rrvert}{\vert}
\newcommand{\llvert}{\vert}
\newtheorem{theorem}{Theorem}[section]
\newremark{rem}{Remark}[section]
\newtheorem{cor}{Corollary}[section]
\newtheorem{lem}{Lemma}[section]
\newtheorem{proposition}{Proposition}[section]
\newproclaim{con}{Condition}[section]

\newcommand{\bbA}{\mathbf{A}}
\newcommand{\bbB}{\mathbf{B}}
\newcommand{\bbD}{\mathbf{D}}
\newcommand{\bbe}{\mathbf{e}}
\newcommand{\bbH}{\mathbf{H}}
\newcommand{\bbI}{\mathbf{I}}
\newcommand{\bbM}{\mathbf{M}}
\newcommand{\bbq}{\mathbf{q}}
\newcommand{\bbr}{\mathbf{r}}
\newcommand{\bbS}{\mathbf{S}}
\newcommand{\bbu}{\mathbf{u}}
\newcommand{\bbW}{\mathbf{W}}
\newcommand{\bbX}{\mathbf{X}}
\newcommand{\bbY}{\mathbf{Y}}

\def\Var{\operatorname{Var}}

\newcommand{\tr}{\operatorname{tr}}
\newcommand{\cov}{\operatorname{cov}}
\newcommand{\diag}{\operatorname{diag}}
\newcommand{\xrightarroww}[1]{\stackrel{#1}{\rule[2.3pt]{0pt}{0.3pt}\hspace*{-1.5pt}\longrightarrow}}
\newcommand{\xrightarrow}[1]{\stackrel{#1}{\rule[2.3pt]{5pt}{0.3pt}\hspace*{-4pt}\longrightarrow}}

\newcommand{\eqref}[1]{(\ref{#1})}
\renewcommand{\epsilon}{\varepsilon}
\renewcommand{\pi}{\uppi}
\newcommand{\dd}{\,\mathrm{d}}

\def\sfrac#1#2{#1/#2}

\def\afrac#1#2{#1/(#2)}

\def\sklfrac#1#2{(#1/#2)}

\newcommand{\diagg}{\mathrm{diag}}
\newcommand{\trr}{\mathrm{tr}}

\makeatother

\begin{document}
\begin{frontmatter}

\title{CLT for linear spectral statistics of normalized sample
covariance matrices with the dimension much larger than the sample size}
\runtitle{CLT for normalized sample covariance matrices}

\begin{aug}
\author[A]{\inits{B.}\fnms{Binbin} \snm{Chen}\thanksref{e1}\ead[label=e1,mark]{CHEN0635@e.ntu.edu.sg}} \and
\author[A]{\inits{G.}\fnms{Guangming} \snm{Pan}\corref{}\thanksref{e2}\ead[label=e2,mark]{gmpan@ntu.edu.sg}}
\address[A]{Division of Mathematical Sciences, School of Physical and
Mathematical Sciences, Nanyang Technological University, Singapore.
\printead{e1,e2}}
\end{aug}

\received{\smonth{1} \syear{2012}}
\revised{\smonth{1} \syear{2014}}

%
\begin{abstract}
Let $\mathbf{A}=\frac{1}{\sqrt{np}}(\mathbf{X}^{T}\mathbf{X}-p\mathbf{I}_n)$ where $\mathbf{X}
$ is a $p\times n$ matrix, consisting of independent and identically
distributed (i.i.d.) real random variables $X_{ij}$ with mean zero
and variance one. When $p/n\to\infty$, under fourth moment conditions
a central limit theorem (CLT) for linear spectral statistics (LSS) of
$\mathbf{A}$ defined by the eigenvalues is established. We also explore its
applications in testing whether a population covariance matrix is an
identity matrix.
\end{abstract}

%

\end{frontmatter}

\section{Introduction}
\label{intoduction}

The last few decades have seen explosive growth in data analysis, due
to the rapid development of modern information technology. We are now
in a setting where many very important data analysis problems are
high-dimensional. In many scientific areas, the data dimension can even
be much larger than the sample size. For example, in micro-array
expression, the number of genes can be tens of thousands or hundreds of
thousands while there are only hundreds of samples. Such kind of data
also arises in genetic, proteomic, functional magnetic resonance
imaging studies and so on (see Chen \textit{et al.} \cite{CZZ10},
Donoho \cite{Donoho00}, Fan and Fan \cite{FF08}).

The main purpose of this paper is to establish a central limit theorem
(CLT) of linear functionals of eigenvalues of the sample covariance
matrix when the dimension $p$ is much larger than the sample size $n$.
Consider the sample covariance matrix
$\bbS=\frac{1}{n}\bbX\bbX^{T}$, where
$\bbX=(X_{ij})_{p\times n}$ and
$X_{ij},i=1,\ldots,p,j=1,\ldots,n$ are i.i.d. real random variables with
mean zero and variance one. As we know, linear functionals of
eigenvalues of $\bbS$ are closely related to its empirical spectral
distribution (ESD) function $F^{\bbS}(x)$. Here for any $n\times n$
Hermitian matrix $\bbM$ with real eigenvalues $\lambda_1,\dots
,\lambda_n$, the empirical spectral distribution of $\bbM$ is defined by
\[
F^{\bbM}=\frac{1}{n}\sum_{j=1}^{n}I(
\lambda_j\leq x),
\]
where $I(\cdot)$ is the indicator function. However, it is
inappropriate to investigate $F^{\bbS}(x)$ when $p/n\to\infty$ since
$\bbS$ has $(p-n)$ zero eigenvalues and hence $F^{\bbS}(x)$ converges
to a degenerate distribution with probability one. Note that the
eigenvalues of $\bbS$ are the same as those of $\frac{1}{n}\bbX
^T\bbX$ except $(p-n)$ zero eigenvalues.
Thus, instead, we turn to the eigenvalues of $\frac{1}{p}\bbX^T\bbX$
and re-normalize it as
%
\begin{equation}
\label{v1} \bbA=\sqrt{\frac{p}{n}} \biggl(\frac{1}{p}
\bbX^T\bbX-\bbI_n \biggr),
\end{equation}
where $\bbI_n$ is the identity matrix of order $n$. 

The first breakthrough regarding the ESD of $\bbA$ was made in
Bai \textit{et al.} \cite{BY88}. They proved that with probability one
\[
F^{\bbA}(x)\rightarrow F(x),
\]
which is the so-called semicircle law with the density
%
\begin{eqnarray}
\label{a5} F'(x)=\lleft\{ %
\begin{array} {l@{\qquad}l}
\dfrac{1}{2\pi}\sqrt{4-x^{2}}, & \mbox{if } |x| \leq2,
\\
0, & \mbox{if } |x| > 2.
\\
\end{array} %
\rright.
\end{eqnarray}
In random matrix theory, $F(x)$ is referred to as the limiting spectral
distribution (LSD) of $\bbA$. For such matrices, Chen and Pan
\cite{CP12} proved
that the largest eigenvalue converges to the right endpoint of the
support of $F(x)$ with probability one. When $X_{11}\sim N(0,1)$,
Karoui \cite{Kar03} reported that the largest eigenvalue
of $\bbX\bbX^T$ after
properly centering and scaling converges in distribution to the
Tracy--Widom law, and Birke and Dette \cite{BD05} established
central limit theorems
for the quadratic function of the eigenvalues of $\bbA$. Recently,
Pan and Gao \cite{PG12} further derived the LSD of a general
form of (\ref{v1}),
which is determined by its Stieltjes transform. Here, the Stieltjes
transform for any distribution function $G$ is given by
\[
m_G(z)= \int\frac{1}{x-z}\dd G(x), \qquad \Im(z)>0,
\]
where $\Im(z)$ represents the imaginary part of $z$.


Gaussian fluctuations in random matrices are investigated by different
authors, starting with Costin and Lebowitz \cite{CL95}. Johansson \cite{Jo.K98} considered an
extended random ensembles whose entries follow a specific class of
densities and established a CLT of the linear spectral statistics
(LSS). Recently, a CLT for LSS of sample covariance matrices is studied
by Bai and Silverstein \cite{BS04} and of Wigner matrices is studied by
Bai and Yao \cite{BY05}.


We introduce some notation before stating our results. Denote the
Stieltjes transform of the semicircle law $F$ by $m(z)$. $\Im(z)$ is
used to denote the imaginary part of a complex number $z$. For any
given square matrix $\bbB$, let $\tr\bbB$ and $\overline{\bbB}$
denote the trace and the\vadjust{\goodbreak} complex conjugate matrix of $\bbB$,
respectively. The norm $\|\bbB\|$ represents the spectral norm of
$\bbB$, that is, $\|\bbB\|=\sqrt{\lambda_{1}(\bbB\overline{\bbB})}$
where $\lambda_{1}(\bbB\overline{\bbB})$ means the maximum
eigenvalue of $\bbB\overline{\bbB}$. The notation $\xrightarroww{d}$
means ``convergence in distribution to''. 
Let $\mathscr{S}$ denote any open region on the real plane including
$[-2,2]$, which is the support of $F(x)$, and $\mathscr{M}$ be the set
of functions which are analytic on $\mathscr{S}$. For any $f\in
\mathscr{M}$, define
%
\begin{equation}
\label{z22} G_n(f)\triangleq n\int_{-\infty}^{+\infty}f(x)\,\mathrm{d}
\bigl(F^{\bbA
}(x)-F(x) \bigr)-\frac{n}{2\pi \mathrm{i}}\oint_{|m|=\rho
}f
\bigl(-m-m^{-1}\bigr)\mathcal{X}_n(m)\frac{1-m^2}{m^2}\dd
m,
\end{equation}
where
%
\begin{eqnarray}
\label{y4}
 \mathcal{X}_n(m) &\triangleq&
 \frac{-\mathcal{B}+\sqrt{\mathcal
{B}^2-4\mathcal{A}\mathcal{C}}}{2\mathcal{A}},\qquad   \mathcal{A} = m-\sqrt{ \frac{n}{p}}\bigl(1+m^2
\bigr),\nonumber
\\[-8pt]\\[-8pt]
\mathcal{B} &=& m^2-1-\sqrt{ \frac{n}{p}}m\bigl(1+2m^2
\bigr), \qquad \mathcal{C} =  \frac{m^3}{n} \biggl( \frac{m^2}{1-m^2}+
\nu_4-2 \biggr)-\sqrt{ \frac{n}{p}}m^4,\nonumber
\end{eqnarray}
$\nu_4 = EX_{11}^{4}$ and $\sqrt{\mathcal{B}^2-4\mathcal{A}\mathcal
{C}}$ is a complex number whose imaginary part has the same sign as
that of $\mathcal{B}$. The integral's contour is taken as $|m|=\rho$
with $\rho<1$.

Let $\{T_k\}$ be the family of Chebyshev polynomials, which is defined
as $T_0(x)=1,T_1(x)=x$ and $T_{k+1}(x)=2xT_k(x)-T_{k-1}(x)$. To give an
alternative way of calculating the asymptotic covariance of $X(f)$ in
Theorem~\ref{thm2} below, for any $f\in\mathscr{M}$ and any integer
$k>0$, we define
\begin{eqnarray*}
\Psi_k(f) & \triangleq&\frac{1}{2\pi}\int
_{-\pi}^{\pi}f(2\cos \theta)\mathrm{e}^{\mathrm{i}k\theta}
\dd\theta
\\
& =&\frac{1}{2\pi}\int_{-\pi}^{\pi}f(2\cos\theta)
\cos k\theta \dd\theta=\frac{1}{\pi}\int_{-1}^{1}f(2x)T_k(x)
\frac{1}{\sqrt{1-x^2}}\dd x. %
\end{eqnarray*}
The main result is formulated in the following.

%
\begin{theorem}\label{thm2} Suppose that
\begin{enumerate}[(b1)]
\item[(a)] $\bbX=(X_{ij})_{p\times n}$ where
$\{X_{ij}\dvtx i=1,2,\ldots,p;j=1,2,\ldots,n\}$ are i.i.d. real random
variables with $EX_{11}=0,EX_{11}^{2}=1 $ and $\nu
_4=EX_{11}^{4}<\infty$.

\item[(b1)] $n/p\to0$ as $n\to\infty$.
\end{enumerate}
Then, for any $f_1,\dots,f_k\in\mathscr{M}$, the finite-dimensional
random vector $ (G_n(f_1),\dots,G_n(f_k) )$ 
converges weakly to a Gaussian vector $ (Y(f_1),\dots,Y(f_k)
)$ with mean function $EY(f)=0$ and covariance function
%
\begin{eqnarray}
\label{ad16}\cov\bigl(Y(f_1),Y(f_2)\bigr)&=&(\nu_4-3)
\Psi_1(f_1)\Psi_1(f_2)+2\sum
_{k=1}^{\infty}k\Psi_k(f_1)
\Psi_k(f_2)
\\
\label{x6}&=&\frac{1}{4\pi^2}\int_{-2}^2\int
_{-2}^2f_1'(x)f_2'(y)H(x,y)\dd x\dd y,
\end{eqnarray}
where
\[
H(x,y)=(\nu_4-3)\sqrt{4-x^2}\sqrt{4-y^2}+2
\log{ \biggl(\frac
{4-xy+\sqrt{(4-x^2)(4-y^2)}}{4-xy-\sqrt{(4-x^2)(4-y^2)}} \biggr)}.
\]
\end{theorem}

%


%
\begin{rem}
Note that $\mathcal{X}_n(m)$ in \eqref{z22} and $\underline{\mathcal
{X}}_n(m) \triangleq\frac{-\mathcal{B}-\sqrt{\mathcal
{B}^2-4\mathcal{A}\mathcal{C}}}{2\mathcal{A}}$ are the two roots of
the equation $\mathcal{A}x^2+\mathcal{B}x+\mathcal{C}=0$. Since
$n/p\to0$, an easy calculation shows $\mathcal{X}_n(m)=\mathrm{o}(1)$ and
$\underline{\mathcal{X}}_n(m)=\frac{1-m^2}{m}+\mathrm{o}(1)$. Hence in
practice, one may implement the mean correction in \eqref{z22} by taking
\[
\mathcal{X}_n(m) = \min \biggl\{\biggl\llvert \frac{-\mathcal{B}+\sqrt {\mathcal{B}^2-4\mathcal{A}\mathcal{C}}}{2\mathcal{A}}
\biggr\rrvert , \biggl\llvert \frac{-\mathcal{B}-\sqrt{\mathcal{B}^2-4\mathcal
{A}\mathcal{C}}}{2\mathcal{A}}\biggr\rrvert \biggr\},
\]
and $m=\rho \mathrm{e}^{\mathrm{i}\theta}$ with $\theta\in[-2\pi,2\pi]$ and $0< \rho
< 1$.
\end{rem}

The mean correction term, the last term in \eqref{z22}, can be
simplified when $n^3/p=\mathrm{O}(1)$. Indeed,
if $n^3/p=\mathrm{O}(1)$, we have $4\mathcal{A}\mathcal{C}= \mathrm{o}(1)$, $\mathcal
{B}=m^2-1$. By \eqref{y4},
\begin{eqnarray*}
n\mathcal{X}_n(m) &=& n\cdot\frac{-\mathcal{B}+\sqrt{\mathcal
{B}^2-4\mathcal{A}\mathcal{C}}}{2\mathcal{A}} =
\frac
{-2nC}{\mathcal{B}+\sqrt{\mathcal{B}^2-4\mathcal{A}\mathcal{C}}}
\\
& =& \frac{m^3}{1-m^2} \biggl(\frac{m^2}{1-m^2} - \nu_4-2 \biggr) +
\sqrt {\frac{n^3}{p}}\frac{m^4}{1-m^2} + \mathrm{o}(1). %
\end{eqnarray*}
Hence, by using the same calculation as that in Section~5.1 of
Bai and Yao \cite{BY05}, we have
%
\begin{eqnarray}
\label{y3} %
&&-\frac{n}{2\pi \mathrm{i}}\oint_{|m|=\rho}f
\bigl(-m-m^{-1}\bigr)\mathcal{X}_n(m)\frac
{1-m^2}{m^2}\dd
m\nonumber
\\
&&\quad =-\frac{1}{2\pi \mathrm{i}}\oint_{|m|=\rho}f\bigl(-m-m^{-1}\bigr)m
\biggl[\frac
{m^2}{1-m^2} - \nu_4-2+\sqrt{\frac{n^3}{p}}m \biggr]
\dd m + \mathrm{o}(1)\nonumber
\\
&&\quad =-\frac{1}{4}\bigl(f(2)+f(-2)\bigr) - \frac{1}{\pi}\int
_{-1}^1f(2x) \biggl[2(\nu_4-3)x^2-
\biggl(\nu_4-\frac{5}{2}\biggr) \biggr]\frac{1}{\sqrt{1-x^2}}\dd x\qquad
\\
&&\qquad {} - \frac{1}{\pi}\sqrt{\frac{n^3}{p}}\int_{-1}^1f(2x)
\frac
{4x^3-3x}{\sqrt{1-x^2}}\dd x\nonumber
\\
&&\quad =- \biggl[ \frac{1}{4}\bigl(f(2)+f(-2)\bigr) - \frac{1}{2}
\Psi_0(f) + (\nu _4-3)\Psi_2(f) \biggr] -
\sqrt{\frac{n^3}{p}}\Psi_3(f) +\mathrm{o}(1). \nonumber%
\end{eqnarray}
%
Define
%
\begin{equation}
\label{z2} Q_n(f)\triangleq n\int_{-\infty}^{+\infty}f(x)\dd
\bigl(F^{\bbA
}(x)-F(x) \bigr)- \sqrt{\frac{n^3}{p}}
\Psi_3(f).
\end{equation}
Under the condition $n^3/p = \mathrm{O}(1)$, we then give a simple and explicit
expression of the mean correction term of \eqref{z22} in the following
corollary.

%
\begin{cor}\label{cor1} Suppose that
\begin{enumerate}[(b2)]
\item[(a)] $\bbX=(X_{ij})_{p\times n}$ where
$\{X_{ij}\dvtx i=1,2,\ldots,p;j=1,2,\ldots,n\}$ are i.i.d. real random
variables with $EX_{11}=0,EX_{11}^{2}=1 $ and $\nu
_4=EX_{11}^{4}<\infty$.

\item[(b2)] $n^3/p = \mathrm{O}(1)$ as $n\to\infty$.
\end{enumerate}
Then, for any $f_1,\dots,f_k\in\mathscr{M}$, the finite-dimensional
random vector $ (Q_n(f_1),\dots,Q_n(f_k) )$ 
converges weakly to a Gaussian vector $ (X(f_1),\dots,X(f_k)
)$ with mean function
%
\begin{equation}
\label{x1} %
EX(f)=\tfrac{1}{4}
\bigl(f(2)+f(-2)\bigr) - \tfrac{1}{2}\Psi_0(f) + (\nu
_4-3)\Psi_2(f) %
\end{equation}
and covariance function $\cov(X(f),X(g))$ being the same as that given
in \eqref{ad16} and \eqref{x6}.
\end{cor}

%
%

%
\begin{rem}\label{rmk3}
The result of Bai \textit{et al.} \cite{Baietal} suggests that, for large
$p$ and $n$ with
$p/n\to\infty$, the matrix $\sqrt{n}\bbA$ is close to a $n\times n$
Wigner matrix although its entries are not independent but weakly
dependent. It is then reasonable to conjecture that the CLT for the LSS
of $\bbA$ resembles that of a Wigner matrix described in Bai and Yao \cite{BY05}.
More precisely, by writing $\bbA=\frac{1}{\sqrt{n}}(w_{ij})$, where
$w_{ii}=(\mathbf{s}_i^T\mathbf{s}_i-p)/\sqrt{p}$, $w_{ij}=\mathbf
{s}_i^T\mathbf{s}_j/\sqrt{p}$ for $i\neq j$ and $\mathbf{s}_j$ is
the $j$th column of
$\bbX$, we have
\[
\Var(w_{11})=\nu_4-1,\qquad  \Var(w_{12})=1, \qquad E
\bigl(w_{12}^2-1\bigr)^2=\frac
{1}{p}
\bigl(\nu_4^2-3\bigr).
\]
%
Then, \eqref{x1}, \eqref{ad16} and \eqref{x6} are consistent with
(1.4), (1.5) and (1.6) of Bai and Yao \cite{BY05},
respectively, by taking their
parameters as $\sigma^2=\nu_4-1$, $\kappa=2$ (the real variable
case) and $\beta=0$.

However, we remark that the mean correction term of $Q_n(f)$, the last
term of (\ref{z2}), cannot be speculated from the result of
Bai and Yao \cite{BY05}. Note that this correction term will
vanish in the case of the
function $f$ to be even or $n^3 /p\rightarrow0$.
By the definition of $\Psi_k(f)$, one may verify that
\begin{eqnarray*}
\Psi_3(f)&=& \frac{1}{\pi}\sqrt{\frac{n^3}{p}}\int
_{-1}^1f(2x)\frac
{4x^3-3x}{\sqrt{1-x^2}}\dd x,
\\
- \frac{1}{2}\Psi_0(f) + (\nu_4-3)
\Psi_2(f)&=&\frac{1}{\pi}\int_{-1}^1f(2x)
\biggl[2(\nu_4-3)x^2-\biggl(\nu_4-
\frac{5}{2}\biggr) \biggr]\frac
{1}{\sqrt{1-x^2}}\dd x.
\end{eqnarray*}
\end{rem}

\begin{rem}\label{rmk4}
If we interchange the roles of $p$ and $n$, Birke and Dette \cite
{BD05} established
the CLT for $Q_n(f)$ in their Theorem~3.4 when $f=x^2$ and $X_{ij}\sim
N(0,1)$. We below show that our Corollary~\ref{cor1} can recover their
result. First, since $f=x^2$ is an even function, it implies that the
last term of (\ref{z2}) is exactly zero. Therefore, the mean in
Theorem~3.4 of Birke and Dette \cite{BD05} is the same as (\ref
{x1}), which equals
one. Second, the variance in Theorem~3.4 of Birke and Dette \cite
{BD05} is also
consistent with (\ref{ad16}). In fact, the variance of Birke and Dette \cite{BD05}
equals $4$ when taking their parameter $y=0$. On the other hand,\vadjust{\goodbreak} since
$X_{ij}\sim N(0,1)$, we have $\nu_4=3$ and the first term of (\ref
{ad16}) is zero. Furthermore, by a direct evaluation, we obtain that
\begin{eqnarray*}
\Psi_1(f) &=& \frac{1}{2\pi}\int_{-\pi}^{\pi}4
\cos^3\theta \dd\theta=\frac{1}{2\pi}\int_{-\pi}^{\pi}
(\cos3\theta+3\cos \theta )\dd\theta=0,
\\
\Psi_2(f) &=& \frac{1}{2\pi}\int_{-\pi}^{\pi}4
\cos^2\theta\cos 2\theta \dd\theta=\frac{1}{2\pi}\int
_{-\pi}^{\pi} (\cos4\theta +1+2\cos2\theta )\dd
\theta=1,
\\
\Psi_k(f)&=&\frac{1}{2\pi}\int_{-\pi}^{\pi}4
\cos^2\theta\cos k\theta \dd\theta=\frac{1}{2\pi}\int
_{-\pi}^{\pi}2 (\cos 2\theta+1 )\cos k\theta \dd\theta
\\
&=&\frac{1}{2\pi}\int_{-\pi}^{\pi} \bigl(\cos(k-2)
\theta+\cos (k+2)\theta+2\cos k\theta \bigr)\dd\theta=0,\qquad  \mbox{for } k\geq3.
\end{eqnarray*}
It implies that $\cov(X(x^2),X(x^2))=4$, which equals the variance of
Birke and Dette \cite{BD05}.
\end{rem}

%



The main contribution of this paper is summarized as follows. We have
established the central limit theorems of linear spectral statistics of
the eigenvalues of the normalized sample covariance matrices when both
the dimension and the sample size go to infinity with the dimension
dominating the sample size (for the case $p/n\to\infty$). Theorem~\ref{thm2} and Corollary~\ref{cor1} are both applicable to the data
with the dimension dominating the sample size while Corollary~\ref
{cor1} provides a simplified correction term (hence, CLT) in the
ultrahigh dimension cases ($n^3/p=\mathrm{O}(1)$). Such an asymptotic theory
complements the results of Bai and Silverstein \cite{BS04} and Pan \cite{Pan} for the case
$p/n\to c\in
(0,\infty)$ and Bai and Yao \cite{BY05} for Wigner matrix.


This paper is organized as follows. Section~\ref{Nur} provides a
calibration of the mean correction term in \eqref{z22}, runs
simulations to check the accuracy of the calibrated CLTs in Theorem~\ref{thm2}, and considers a statistical
application of Theorem~\ref{thm2} and a real data analysis. Section~\ref{TandS} gives the strategy of proving Theorem~\ref{thm2} and two
intermediate results, Propositions \ref{pro1} and \ref{pro3}, and
truncation steps of the underlying random variables are given as well.
Some preliminary results are given in Section~\ref{lemmas}. Sections~\ref{M1} and \ref{M2} are devoted to the proof of Proposition~\ref
{pro1}. We present the proof of Proposition~\ref{pro3} in Section~\ref{proofPro3}. Section~\ref{calculation} derives mean and covariance in
Theorem~\ref{thm2}.


\section{Calibration, application and empirical studies}\label{Nur}

Section~\ref{Calib} considers a calibration to the mean correction
term of \eqref{z22}. A statistical application is performed in
Section~\ref{test} and the empirical studies are carried out in
Section~\ref{empirical}.

\subsection{Calibration of the mean correction term in \texorpdfstring{\protect\eqref{z22}}{(1.3)}}
\label{Calib}

Theorem~\ref{thm2} provides a CLT for $G_n(f)$ under the general
framework $p/n\to\infty$, which only requires zero mean, unit
variance and the bounded fourth moment. However, the simulation results
show that the asymptotic distributions of $G_n(f)$, especially the
asymptotic means, are sensitive to the skewness and the kurtosis of the
random variables for some particular functions~$f$, for example, $f(x)=\frac
{1}{2}x(x^2-3)$. This phenomenon is caused by the slow convergence rate
of $EG_n(f)$ to zero, which is illustrated as follows. Suppose that
$EX_{11}^8<\infty$. We then have $|EG_n(f)|=\mathrm{O}(\sqrt{n/p})+\mathrm{O}(1/\sqrt {n})$ by the arguments in Section~\ref{M2}. Also, the remaining terms
(see \eqref{remain} below) have a coefficient $(\nu_4-1)\sqrt{n/p}$
which converges to zero theoretically since $p/n\rightarrow\infty$.
However, if $n = 100, p = n^2$ and the variables $X_{ij}$ are from
central $\exp(1)$ then $(\nu_4-1)\sqrt{n/p}$ could be as big as $0.8$.


In view of this, we will regain such terms and give a calibration for
the mean correction term in \eqref{z22}. From Section~\ref{M2}, we
observe that the convergence rate of $|EG_n(f)|$ relies on the rate of
$|nE\omega_n - m^3(z)(m'(z)+\nu_4-2)|$ in Lemma~\ref{lem4}. By the
arguments in Section~\ref{M2}, only $S_{22}$ below \eqref{s23} has
the coefficient $(\nu_4-1)\sqrt{n/p}$. A simply calculation implies
that
\begin{equation}
\label{remain} %
S_{22} = -2(\nu_4-1)
\sqrt{n/p}m(z)+\mathrm{o}(1). %
\end{equation}
%
Hence, the limit of $nE\omega_n$ is calibrated as
%
\begin{equation}
\label{o1} nE\omega_n = m^3(z) \bigl[
\nu_4-2+m'(z)-2(\nu_4-1)\sqrt{n/p}m(z)
\bigr] + \mathrm{o}(1).
\end{equation}
We then calibrate $G_n(f)$ as
%
\begin{eqnarray}
\label{o2} G_n^{\mathrm{Calib}}(f)&\triangleq& n\int
_{-\infty}^{+\infty}f(x)\,\mathrm{d} \bigl(F^{\bbA}(x)-F(x)
\bigr)\nonumber\\[-8pt]\\[-8pt]
&&{}-\frac{n}{2\pi \mathrm{i}}\oint_{|m|=\rho
}f\bigl(-m-m^{-1}\bigr)
\mathcal{X}_n^{\mathrm{Calib}}(m)\frac{1-m^2}{m^2}\dd m,\nonumber
\end{eqnarray}
where, via \eqref{o1},
%
\begin{eqnarray}
\label{Xcal} \mathcal{X}_n^{\mathrm{Calib}}(m) &\triangleq&
\frac{-\mathcal{B}+\sqrt {\mathcal{B}^2-4\mathcal{A}\mathcal{C}^{\mathrm{Calib}}}}{2\mathcal{A}},\nonumber\\[-8pt]\\[-8pt]
 \mathcal{C}^{\mathrm{Calib}} &=& \frac{m^3}{n} \biggl[
\nu_4-2+\frac
{m^2}{1-m^2}-2(\nu_4-1)m\sqrt{n/p} \biggr]-
\sqrt{\frac{n}{p}}m^4,\nonumber
\end{eqnarray}
$\mathcal{A},\mathcal{B}$ are defined in \eqref{y4} and $\sqrt {\mathcal{B}^2-4\mathcal{A}\mathcal{C}^{\mathrm{Calib}}}$ is a complex number
whose imaginary part has the same sign as that of $\mathcal{B}$.
Theorem~\ref{thm2} still holds if we replace $G_n(f)$ with $G_n^{\mathrm{Calib}}(f)$.


We next perform a simulation study to check the accuracy of the CLT in
Theorem~\ref{thm2} with $G_n(f)$ replaced by the calibrated expression
$G_n^{\mathrm{Calib}}(f)$ in \eqref{o2}. Two combinations of $(p,n)$,
$p=n^2,n^{2.5}$, and the test function $f(x)=\frac{1}{2}x(x^2-3)$ are
considered in the simulations, as suggested by one of the referees. To
inspect the impact of the skewness and the kurtosis of the variables,
we use three types of random variables, $N(0,1)$, central $\exp(1)$ and
central $t(6)$. The skewnesses of these variables are 0, 2 and 0 while
the fourth moments of these variables are 3, 9 and 6, respectively. The
empirical means and empirical standard deviations of
$G_n^{\mathrm{Calib}}(f)/(\Var(Y(f)))^{1/2}$ from 1000 independent replications
are shown in Table~\ref{table1}.
\begin{table}
\tablewidth=\textwidth
\tabcolsep=0pt
\caption{Empirical means of $G_n^{\mathrm{Calib}}(f)/(\Var(Y(f)))^{1/2}$ (cf.
\protect\eqref{o2}) for the function $f(x) = {1\over2}x(x^2-3)$ with the
corresponding standard deviations in the parentheses}\label{table1}
\begin{tabular*}{\textwidth}{@{\extracolsep{\fill}}lllll@{}}
\hline
 $n$ & 50 &100 &150 &200\\
\hline
&\multicolumn{4}{l}{$p=n^2$}\\
$N(0,1)$ &$-0.314\ (1.227)$ &$-0.221\ (1.038)$ &$-0.188\ (1.051)$
&$-0.093\ (0.940)$\\
$\exp(1)$ &$-0.088\ (2.476)$ &$-0.079\ (1.447)$ &$-0.140\ (1.400)$
&$-0.161\ (1.154)$\\
$t(6)$ &$-0.084\ (2.813)$ &$-0.077\ (1.541)$ &$-0.095\ (1.246)$
&$-0.0897\ (1.104)$\\[3pt]
& \multicolumn{4}{l}{$p=n^{2.5}$} \\
$N(0,1)$ &$-0.068\ (1.049)$ &$-0.053\ (1.077)$ &$-0.0476\ (0.944)$
&$-0.016\ (1.045)$\\
$\exp(1)$ &$-0.049\ (1.879)$ &$-0.029\ (1.390)$ &$-0.046\ (1.162)$
&$-0.045\ (1.156)$\\
$t(6)$ &$-0.075\ (1.693)$ &$ \hphantom{-}0.050\ (1.252)$ &$-0.044\ (1.145)$
&$-0.027\ (1.044)$\\
\hline
\end{tabular*}
\end{table}
%
\begin{figure}[b]

\includegraphics{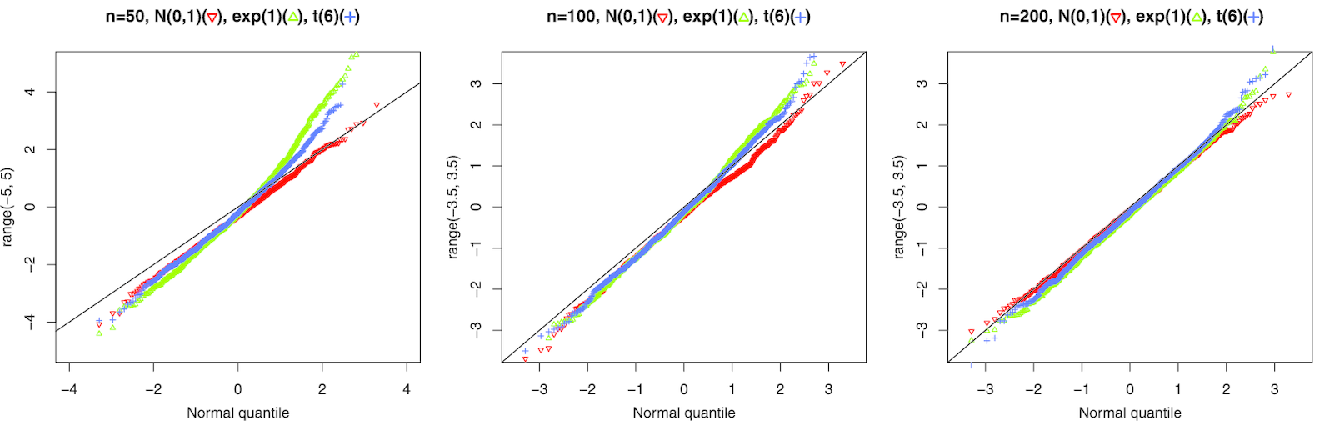}

\caption{The Q--Q plots of the standard Gaussian distribution versus
$G_n^{\mathrm{Calib}}(f)/(\Var(Y(f)))^{1/2}$ based on the sample generating from
$N(0,1)$ ($\triangledown$), standardized $\exp(1)$ ($\triangle$) and standardized $t(6)$
(+) with the sample sizes $n=50,100,200$ from left to right and the
dimension $p=n^2$.}\label{fig1}
\end{figure}
%
\begin{figure}

\includegraphics{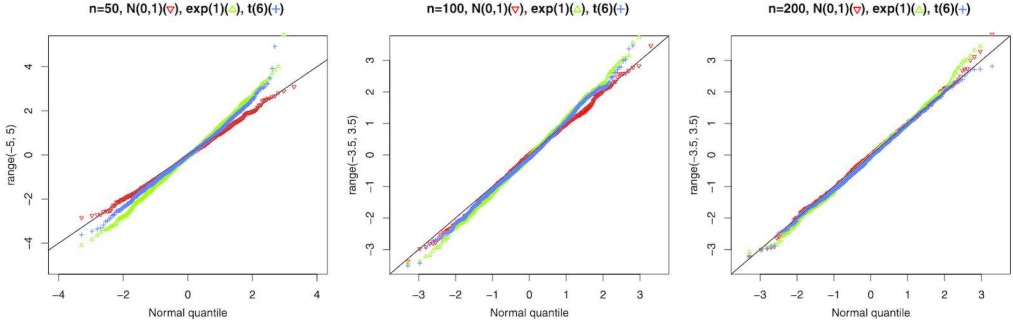}

\caption{The Q--Q plots of the standard Gaussian distribution versus
$G_n^{\mathrm{Calib}}(f)/(\Var(Y(f)))^{1/2}$ based on the sample generating from
$N(0,1)$ ($\triangledown$), standardized $\exp(1)$ ($\triangle$) and standardized $t(6)$
(+) with the sample sizes $n=50,100,200$ from left to right and the
dimension $p=n^{2.5}$.}\label{fig2}
\end{figure}

It is observed from Table~\ref{table1} that both the empirical means and standard
deviations for $N(0,1)$ random variables are very accurate. The
empirical means for central $\exp(1)$ and central $t(6)$ also show their
good accuracy. We note that the standard deviations for central
$\exp(1)$ and central $t(6)$ random variables are not good when $n$ is
small (e.g., $n=50$). But it gradually tends to 1 as the sample size
$n$ increases.

Q--Q plots are employed to illustrate the accuracy of the normal
approximation in Figures \ref{fig1} and \ref{fig2} corresponding to the scenarios $p=n^2$
and $p=n^{2.5}$, respectively. In each figure, Q--Q plots from left to
right correspond to $n=50,100,150,200$, respectively with random
variables generated from $N(0,1)$ ($\triangledown$), central $\exp(1)$ ($\triangle$) and
central $t(6)$ (+). We observe the same phenomenon that the normal
approximation is very accurate for normal variables while the
approximation is gradually better when $n$ increases for central
$\exp(1)$ and $t(6)$ variables.

\subsection{Application of CLTs to hypothesis test}\label{test}
This subsection is to consider an application of Theorem~\ref{thm2}
which is about hypothesis testing for the covariance matrix.
Suppose that $\mathbf{y}=\Gamma\mathbf{s}$ is a $p$-dimensional
vector where
$\Gamma$ is a $p\times p$ matrix with positive eigenvalues and the
entries of $\mathbf{s}$ are i.i.d. random variables with mean zero and
variance one. Hence, the covariance matrix of $\mathbf{y}$ is $\Sigma
=\Gamma\Gamma^T$. Suppose that one wishes to test the hypothesis
%
\begin{equation}
\label{k1} H_0\dvtx  \Sigma=\bbI_p,\qquad  H_1\dvtx
\Sigma\neq\bbI_p.
\end{equation}
%

Based on the i.i.d. samples $\mathbf{y}_1,\dots,\mathbf{y}_n$ (from
$\mathbf{y}$), many authors have considered (\ref{k1}) in terms of the
relationship of $p$ and $n$. For example, John \cite
{John71} and Nagao \cite{Nagao73} considered the fixed-dimensional case; Ledoit and Wolf \cite{LW02}, Fisher
\textit{et al.} \cite{Fisher10} and Bai \textit{et al.}
\cite{Baietal} studied the case of $\frac
{p}{n}\to c\in(0,\infty)$; Srivastava \cite{Sva05},
Srivastava, Kollo and von Rosen \cite{Sva11}, Fisher \cite
{Fisher12} and Chen \textit{et al.} \cite{CZZ10} proposed the testing
statistics which can
accommodate large $p$ and small $n$.

We are interesting in testing (\ref{k1}) in the setting of $\frac
{p}{n}\to\infty$. 
As in Ledoit and Wolf \cite{LW02} and Birke and Dette \cite
{BD05}, we set $f=x^2$. 
We then propose the following test statistic for the hypothesis of
(\ref{k1}):
%
\begin{eqnarray}
\label{o4*} L_n &=& \frac{1}{2} \biggl[n \biggl(\int
x^2 \dd F^{\bbB}(x) - \int x^2 \dd F(x) \biggr)\nonumber\\[-8pt]\\[-8pt]
&&\hphantom{\frac{1}{2} \biggl[}{}-
\biggl(\frac{n}{2\pi \mathrm{i}}\oint_{|m|=\rho
}\bigl(m+m^{-1}
\bigr)^2\mathcal{X}_n^{\mathrm{Calib}}(m)\frac{1-m^2}{m^2}
\dd m \biggr) \biggr],\nonumber
\end{eqnarray}
where $\mathcal{X}_n^{\mathrm{Calib}}(m)$ is given in \eqref{Xcal} and $\bbB
=\sqrt{\frac{p}{n}} (\frac{1}{p}\bbY^T\bbY-\bbI_n )$ is
the normalized sample covariance matrix with $\bbY=(\mathbf
{y}_1,\dots
,\mathbf{y}_n)$.
The asymptotic mean and variance of $L_n$ are $0$ and $1$,
respectively, see Theorem~\ref{thm2} or Remark~\ref{rmk4} for
details. Since there is no close form for the mean correction term in
\eqref{o4*}, we use \emph{Matlab} to calculate this correction term.
It shows that as $n/p\rightarrow0$,
\[
\frac{n}{2\pi \mathrm{i}}\oint_{|m|=\rho}\bigl(m+m^{-1}
\bigr)^2\mathcal {X}_n^{\mathrm{Calib}}(m)\frac{1-m^2}{m^2}
\dd m = \nu_4-2.
\]
We also note the fact that
\[
E \biggl[n\int x^2 \,\mathrm{d}\bigl(F^{\bbB}(x)-F(x)\bigr) \biggr] = E
\bigl[\tr\bbB\bbB^T -n \bigr] = \nu_4-2.
\]
Thus, we use the following test statistic in the simulations:
%
\begin{equation}
\label{o4} L_n = \frac{1}{2} \biggl[n \biggl(\int
x^2 \dd F^{\bbB}(x) - \int x^2 \dd F(x) \biggr)-(
\nu_4-2) \biggr] = \frac{1}{2} \bigl(\tr\bbB\bbB^T
-n -(\nu_4-2) \bigr).
\end{equation}
Since $\Gamma^T\Gamma=\bbI_p$ is equivalent to $\Gamma\Gamma
^T=\bbI_p$, under the null hypothesis $H_0$ in (\ref{k1}), we have
%
\begin{equation}
\label{k4} L_n\xrightarroww{d}N(0,1).
\end{equation}
By the law of large numbers, a consistent estimator of $\nu_4$ is
$\widehat{\nu}_4= \frac{1}{np}\sum_{i,j}Y_{ij}^4$ under the null
hypothesis $H_0$. By Slutsky's theorem, \eqref{k4} also holds if we
replace $\nu_4$ of \eqref{o4} with $\widehat{\nu}_4$.\looseness=-1


The numerical performance of the proposed statistic $L_n$ is carried
out by Monte Carlo simulations. Let $Z_{\alpha/2}$ and $Z_{1-\alpha
/2}$, respectively, be the $100\alpha/2\%$ and $100(1-\alpha/2)\%$
quantiles of the asymptotic null distribution of the test statistic
$L_n$. With $T$ replications of the data set simulated under the null
hypothesis, we calculate the empirical size as
\[
\hat{\alpha}=\frac{\{\# L_n^{\mathrm{null}}\leq Z_{\alpha/2}\}+\{\#
L_n^{\mathrm{null}}> Z_{1-\alpha/2}\}}{T},
\]
where $\#$ denotes the number and $L_n^{\mathrm{null}}$ represents the values of
the test statistic $L_n$ based on the data set simulated under the null
hypothesis. The empirical power is calculated by\looseness=-1
\[
\hat{\beta}=\frac{\{\# L_n^{\mathrm{alter}}\leq Z_{\alpha/2}\}+\{\#
L_n^{\mathrm{alter}}> Z_{1-\alpha/2}\}}{T},\looseness=0
\]
where $L_n^{\mathrm{alter}}$ represents the values of the test statistic $L_n$
based on the data set simulated under the alternative hypothesis. In
our simulations, we fix $T=1000$ as the number of replications and set
the nominal significance level $\alpha=5\%$. By asymptotic normality,
we have $Z_{\alpha/2}=-1.96$ and $Z_{1-\alpha/2}=1.96$.
\begin{table}[b]
\tablewidth=\textwidth
\tabcolsep=0pt
\caption{Empirical sizes of CZZ test and $L_n$ at the significant
level $\alpha=5\%$ for normal, gamma, Bernoulli random vectors}\label{table2}
\begin{tabular*}{\textwidth}{@{\extracolsep{\fill}}lllllllll@{}}
\hline
& \multicolumn{4}{l}{CZZ test} & \multicolumn{4}{l}{$L_n$}\\
 & \multicolumn{4}{l}{$n$} & \multicolumn{4}{l}{$n$} \\[-5pt]
 & \multicolumn{4}{l}{\hrulefill} & \multicolumn{4}{l}{\hrulefill} \\
$p$ & 20 &40 &60 &80 & 20 &40 &60 &80\\
\hline
& \multicolumn{8}{l}{Normal random vectors}\\
\hphantom{10\,}600 &0.069 &0.071 &0.052 &0.052 &0.063 &0.077 &0.066 &0.082\\
\hphantom{1\,}1500 &0.057 &0.059 &0.061 &0.059 &0.055 &0.058 &0.058 &0.062\\
\hphantom{1\,}3000 &0.067 &0.068 &0.057 &0.053 &0.048 &0.067 &0.056 &0.052\\
\hphantom{1\,}5500 &0.064 &0.06 &0.067 &0.058 &0.054 &0.055 &0.071 &0.068\\
\hphantom{1\,}8000 &0.071 &0.062 &0.062 &0.054 &0.055 &0.049 &0.06 &0.059\\
10\,000 &0.055 &0.059 &0.063 &0.06 &0.037 &0.058 &0.057 &0.054\\
[3pt]
& \multicolumn{8}{l}{Gamma random vectors}\\
\hphantom{10\,}600 &0.055 &0.073 &0.056 &0.062 &0.103 &0.119 &0.125 &0.123\\
\hphantom{1\,}1500 &0.064 &0.047 &0.059 &0.059 &0.094 &0.072 &0.072 &0.088\\
\hphantom{1\,}3000 &0.069 &0.071 &0.059 &0.052 &0.066 &0.074 &0.071 &0.061\\
\hphantom{1\,}5500 &0.065 &0.069 &0.048 &0.041 &0.077 &0.073 &0.047 &0.045\\
\hphantom{1\,}8000 &0.069 &0.065 &0.07 &0.053 &0.078 &0.075 &0.063 &0.059\\
10\,000 &0.072 &0.06 &0.06 &0.057 &0.078 &0.082 &0.065 &0.06 \\
[3pt]
& \multicolumn{8}{l}{Bernoulli random vectors}\\
\hphantom{10\,}600 &0.078 &0.079 &0.056 &0.037 &0.048 &0.064 &0.046 &0.037\\
\hphantom{1\,}1500 &0.065 &0.050 &0.051 &0.053 &0.039 &0.040 &0.049 &0.050\\
\hphantom{1\,}3000 &0.048 &0.053 &0.058 &0.060 &0.040 &0.052 &0.052 &0.056\\
\hphantom{1\,}5500 &0.059 &0.061 &0.059 &0.042 &0.040 &0.052 &0.060 &0.040\\
\hphantom{1\,}8000 &0.065 &0.074 &0.065 &0.059 &0.046 &0.052 &0.05 &0.051\\
10\,000 &0.07 &0.057 &0.047 &0.048 &0.044 &0.037 &0.038 &0.047\\
\hline
\end{tabular*}
\end{table}
%
\begin{table}[b]
\tablewidth=\textwidth
\tabcolsep=0pt
\caption{Empirical powers of CZZ test and $L_n$ at the significant
level $\alpha=5\%$ for normal random vectors. Two types of population
covariance matrices are considered. In the first case, $\Sigma
_1=\diag(2\times\mathbf{1}_{[\nu p]},\mathbf{1}_{p-[\nu p]})$ for
$\nu=0.08$ and $\nu=0.25$, respectively. In the second case, $\Sigma
_2=\diag(A_1,\diag (\mathbf{1}_{p-[v_2p]}) )$,
where $A_1$ is a $[v_2p]\times[v_2p]$ tridiagonal symmetric matrix
with diagonal elements equal to 1 and elements beside diagonal all
equal to $v_1$ for $v_1=0.5,v_2=0.8$ and $v_1=0.5,v_2=0.4$, respectively}\label{table3}
\begin{tabular*}{\textwidth}{@{\extracolsep{\fill}}lllllllll@{}}
\hline
& \multicolumn{4}{l}{CZZ test} & \multicolumn{4}{l}{$L_n$}\\
 & \multicolumn{4}{l}{$n$} & \multicolumn{4}{l}{$n$} \\[-5pt]
 & \multicolumn{4}{l}{\hrulefill} & \multicolumn{4}{l}{\hrulefill} \\
$p$ & 20 &40 &60 &80 & 20 &40 &60 &80\\
\hline
& \multicolumn{8}{l}{Normal random vectors ($\nu=0.08$)}\\
\hphantom{10\,}600 &0.186 &0.392 &0.648 &0.826 &0.932 &1 &1 &1\\
\hphantom{1\,}1500 &0.179 &0.397 &0.642 &0.822 &0.999 &1 &1 &1\\
\hphantom{1\,}3000 &0.197 &0.374 &0.615 &0.867 &1.000 &1 &1 &1\\
\hphantom{1\,}5500 &0.225 &0.382 &0.615 &0.85 &1 &1 &1 &1\\
\hphantom{1\,}8000 &0.203 &0.391 &0.638 &0.843 &1 &1 &1 &1\\
10\,000 &0.204 &0.381 &0.639 &0.835 &1 &1 &1 &1\\
[3pt]
& \multicolumn{8}{l}{Normal random vectors ($\nu=0.25$)}\\
\hphantom{10\,}600 &0.571 &0.952 &0.997 &1 &1 &1 &1 &1 \\
\hphantom{1\,}1500 &0.585 &0.959 &1.000 &1 &1 &1 &1 &1\\
\hphantom{1\,}3000 &0.594 &0.961 &1.000 &1 &1 &1 &1 &1\\
\hphantom{1\,}5500 &0.617 &0.954 &1 &1 &1 &1 &1 &1\\
\hphantom{1\,}8000 &0.607 &0.957 &0.999 &1 &1 &1 &1 &1\\
10\,000 &0.595 &0.949 &1 &1 &1 &1 &1 &1\\
[3pt]
& \multicolumn{8}{l}{Normal random vectors ($v_1=0.5,v_2=0.8$)}\\
\hphantom{10\,}600 &0.333 &0.874 &0.997 &1 &0.443 &0.493 &0.492 &0.488 \\
\hphantom{1\,}1500 &0.310 &0.901 &0.999 &1 &0.987 &0.997 &0.997 &0.998\\
\hphantom{1\,}3000 &0.348 &0.889 &0.998 &1 &1.000 &1.000 &1.000 &1.000\\
\hphantom{1\,}5500 &0.382 &0.871 &0.998 &1 &1 &1 &1 &1\\
\hphantom{1\,}8000 &0.33 &0.867 &0.998 &1 &1 &1 &1 &1\\
10\,000 &0.359 &0.868 &0.998 &1 &1 &1 &1 &1\\
[3pt]
& \multicolumn{8}{l}{Normal random vectors ($v_1=0.5,v_2=0.4$)}\\
\hphantom{10\,}600 &0.142 &0.364 &0.668 &0.896 &0.078 &0.089 &0.069 &0.102 \\
\hphantom{1\,}1500 &0.131 &0.354 &0.653 &0.890 &0.220 &0.235 &0.230 &0.226\\
\hphantom{1\,}3000 &0.139 &0.361 &0.662 &0.899 &0.635 &0.660 &0.647 &0.684\\
\hphantom{1\,}5500 &0.148 &0.352 &0.645 &0.898 &0.97 &0.979 &0.989 &0.989\\
\hphantom{1\,}8000 &0.152 &0.36 &0.688 &0.905 &0.981 &0.978 &0.986 &0.989\\
10\,000 &0.137 &0.328 &0.674 &0.886 &1 &1 &1 &1\\
\hline
\end{tabular*}
\end{table}
%
\begin{table}
\tablewidth=\textwidth
\tabcolsep=0pt
\caption{Empirical powers of CZZ test and $L_n$ at the significant
level $\alpha=5\%$ for standardized gamma random vectors. Two types of
population covariance matrices are considered. In the first case,
$\Sigma_1=\diag(2\times\mathbf{1}_{[\nu p]},\mathbf{1}_{p-[\nu p]})$
for $\nu=0.08$ and $\nu=0.25$, respectively. In the second case,
$\Sigma_2=\diag(A_1,\diag (\mathbf
{1}_{p-[v_2p]}) )$, where $A_1$ is a $[v_2p]\times[v_2p]$
tridiagonal symmetric matrix with diagonal elements equal to 1 and
elements beside diagonal all equal to $v_1$ for $v_1=0.5,v_2=0.8$ and
$v_1=0.5,v_2=0.4$, respectively}\label{table4}
\begin{tabular*}{\textwidth}{@{\extracolsep{\fill}}lllllllll@{}}
\hline
& \multicolumn{4}{l}{CZZ test} & \multicolumn{4}{l}{$L_n$}\\
 & \multicolumn{4}{l}{$n$} & \multicolumn{4}{l}{$n$} \\[-5pt]
 & \multicolumn{4}{l}{\hrulefill} & \multicolumn{4}{l}{\hrulefill} \\
$p$ & 20 &40 &60 &80 & 20 &40 &60 &80\\
\hline
& \multicolumn{8}{l}{Gamma random vectors ($\nu=0.08$)}\\
\hphantom{10\,}600 &0.331 &0.638 &0.891 &0.982 &0.999 &1 &1 &1\\
\hphantom{1\,}1500 &0.356 &0.636 &0.901 &0.979 &1 &1 &1 &1\\
\hphantom{1\,}3000 &0.197 &0.383 &0.638 &0.823 &1 &1 &1 &1\\
\hphantom{1\,}5500 &0.178 &0.361 &0.658 &0.845 &1 &1 &1 &1\\
\hphantom{1\,}8000 &0.199 &0.399 &0.642 &0.85 &1 &1 &1 &1\\
10\,000 &0.216 &0.353 &0.636 &0.843 &1 &1 &1 &1\\
[3pt]
& \multicolumn{8}{l}{Gamma random vectors ($\nu=0.25$)}\\
\hphantom{10\,}600 &0.621 &0.943 &1.000 &1 &1 &1 &1 &1 \\
\hphantom{1\,}1500 &0.610 &0.946 &0.999 &1 &1 &1 &1 &1\\
\hphantom{1\,}3000 &0.579 &0.946 &0.997 &1 &1 &1 &1 &1\\
\hphantom{1\,}5500 &0.596 &0.957 &0.999 &1 &1 &1 &1 &1\\
\hphantom{1\,}8000 &0.616 &0.962 &0.999 &1 &1 &1 &1 &1 \\
10\,000 &0.614 &0.955 &0.999 &1 &1 &1 &1 &1\\
[3pt]
& \multicolumn{8}{l}{Gamma random vectors ($v_1=0.5,v_2=0.8$)}\\
\hphantom{10\,}600 &0.192 &0.871 &0.998 &0.972 &0.122 &0.413 &0.423 &0.133 \\
\hphantom{1\,}1500 &0.198 &0.883 &0.995 &0.980 &0.440 &0.992 &0.993 &0.433\\
\hphantom{1\,}3000 &0.343 &0.885 &0.995 &1 &1 &1 &1 &1\\
\hphantom{1\,}5500 &0.342 &0.88 &0.996 &1 &1 &1 &1 &1\\
\hphantom{1\,}8000 &0.349 &0.877 &0.998 &1 &1 &1 &1 &1 \\
10\,000 &0.337 &0.879 &0.998 &1 &1 &1 &1 &1\\
[3pt]
& \multicolumn{8}{l}{Gamma random vectors ($v_1=0.5,v_2=0.4$)}\\
\hphantom{10\,}600 &0.117 &0.353 &0.650 &0.780 &0.087 &0.111 &0.114 &0.120 \\
\hphantom{1\,}1500 &0.138 &0.365 &0.661 &0.799 &0.183 &0.215 &0.226 &0.157\\
\hphantom{1\,}3000 &0.129 &0.349 &0.646 &0.89 &0.593 &0.621 &0.627 &0.61\\
\hphantom{1\,}5500 &0.124 &0.335 &0.678 &0.889 &0.945 &0.972 &0.981 &0.986\\
\hphantom{1\,}8000 &0.142 &0.369 &0.668 &0.901 &0.999 &1 &1 &1 \\
10\,000 &0.142 &0.336 &0.668 &1 &1 &1 &1 &1\\
\hline
\end{tabular*}
\end{table}
%
\begin{table}
\tablewidth=\textwidth
\tabcolsep=0pt
\caption{Empirical powers of CZZ test and $L_n$ at the significant
level $\alpha=5\%$ for standardized Bernoulli random vectors. Two
types of population covariance matrices are considered. In the first
case, $\Sigma_1=\diag(2\times\mathbf{1}_{[\nu p]},\mathbf
{1}_{p-[\nu
p]})$ for $\nu=0.08$ and $\nu=0.25$, respectively. In the second
case, $\Sigma_2=\diag(A_1,\diag (\mathbf
{1}_{p-[v_2p]}) )$, where $A_1$ is a $[v_2p]\times[v_2p]$
tridiagonal symmetric matrix with diagonal elements equal to 1 and
elements beside diagonal all equal to $v_1$ for $v_1=0.5,v_2=0.8$ and
$v_1=0.5,v_2=0.4$, respectively}\label{table5}
\begin{tabular*}{\textwidth}{@{\extracolsep{4in minus 4in}}lllllllll@{}}
\hline
& \multicolumn{4}{l}{CZZ test} & \multicolumn{4}{l}{$L_n$}\\
 & \multicolumn{4}{l}{$n$} & \multicolumn{4}{l}{$n$} \\[-5pt]
 & \multicolumn{4}{l}{\hrulefill} & \multicolumn{4}{l}{\hrulefill} \\
$p$ & 20 &40 &60 &80 & 20 &40 &60 &80\\
\hline
& \multicolumn{8}{l}{Bernoulli random vectors ($\nu=0.08$)}\\
\hphantom{10\,}600 &0.216 &0.381 &0.622 &0.849 &0.972 &1 &1 &1\\
\hphantom{1\,}1500 &0.198 &0.401 &0.632 &0.837 &1 &1 &1 &1\\
\hphantom{1\,}3000 &0.203 &0.362 &0.622 &0.823 &1 &1 &1 &1\\
\hphantom{1\,}5500 &0.196 &0.354 &0.627 &0.829 &1 &1 &1 &1\\
\hphantom{1\,}8000 &0.203 &0.373 &0.638 &0.834 &1 &1 &1 &1\\
10\,000 &0.213 &0.397 &0.637 &0.822 &1 &1 &1 &1\\
[3pt]
& \multicolumn{8}{l}{Bernoulli random vectors ($\nu=0.25$)}\\
\hphantom{10\,}600 &0.594 &0.952 &0.998 &1 &1 &1 &1 &1\\
\hphantom{1\,}1500 &0.619 &0.960 &1.000 &1 &1 &1 &1 &1\\
\hphantom{1\,}3000 &0.594 &0.964 &0.999 &1 &1 &1 &1 &1\\
\hphantom{1\,}5500 &0.609 &0.948 &1.000 &1 &1 &1 &1 &1\\
\hphantom{1\,}8000 &0.589 &0.952 &1 &1 &1 &1 &1 &1\\
10\,000 &0.603 &0.957 &0.999 &1 &1 &1 &1 &1\\
[3pt]
& \multicolumn{8}{l}{Bernoulli random vectors ($v_1=0.5,v_2=0.8$)}\\
\hphantom{10\,}600 &0.356 &0.870 &0.996 &1 &0.507 &0.512 &0.526 &0.558\\
\hphantom{1\,}1500 &0.359 &0.892 &0.995 &1 &0.999 &1 &1 &0.999\\
\hphantom{1\,}3000 &0.343 &0.877 &0.998 &1 &1.000 &1 &1 &1.000\\
\hphantom{1\,}5500 &0.355 &0.868 &0.997 &1 &1.000 &1.000 &1.000 &1.000\\
\hphantom{1\,}8000 &0.332 &0.873 &0.997 &1 &1 &1 &1 &1\\
10\,000 &0.353 &0.872 &1 &1 &1 &1 &1 &1\\
[3pt]
& \multicolumn{8}{l}{Bernoulli random vectors ($v_1=0.5,v_2=0.4$)}\\
\hphantom{10\,}600 &0.153 &0.348 &0.643 &0.901 &0.092 &0.086 &0.079 &0.085\\
\hphantom{1\,}1500 &0.154 &0.372 &0.643 &0.878 &0.239 &0.255 &0.235 &0.241\\
\hphantom{1\,}3000 &0.141 &0.339 &0.649 &0.882 &0.682 &0.680 &0.680 &0.674\\
\hphantom{1\,}5500 &0.156 &0.343 &0.656 &0.893 &0.997 &0.994 &0.994 &0.994\\
\hphantom{1\,}8000 &0.144 &0.353 &0.664 &0.904 &1 &1 &1 &1\\
10\,000 &0.139 &0.356 &0.685 &0.889 &1 &1 &1 &1\\
\hline
\end{tabular*}
\end{table}


Our proposed test is intended for the situation ``large $p$, small
$n$''. To inspect the impact caused by the sample size and/or the
dimension, we set
\begin{eqnarray*}
n&=&20,40,60,80,
\\
p&=&600,1500,3000,5500,8000,10\,000.
\end{eqnarray*}
The entries of $\mathbf{s}$ are generated from three types of
distributions, Gaussian distribution, standardized Gamma(4,\ 0.5) and
Bernoulli distribution with $P(X_{ij}=\pm1)=0.5$.\vadjust{\goodbreak} 

The following two types of covariance matrices are considered in the
simulations to investigate the empirical power of the test.
\begin{enumerate}[2.]
\item[1.] (Diagonal covariance.) $\Sigma=\diag(\sqrt{2}\mathbf
{1}_{[\nu p]},\mathbf{1}_{1-[\nu p]})$, where $\nu=0.08$ or $\nu
=0.25$, $[a]$ denotes the largest integer that is not greater than $a$.
\item[2.] (Banded covariance.) $\Sigma=\diag (A_1,\diag(\mathbf
{1}_{p-[v_2p]}) )$, where $A_1$ is a $[v_2p]\times
[v_2p]$ tridiagonal symmetric matrix with the diagonal elements being
equal to 1 and elements below and above the diagonal all being equal to $v_1$.
\end{enumerate}

Since the test in Chen \textit{et al.} \cite{CZZ10} accommodates a
wider class of variates
and has less restrictions on the ratio $p/n$, we below compare
performance of our test with that of Chen \textit{et al.} \cite{CZZ10}.
To simplify the
notation, denote their test by the CZZ test. Table~\ref{table2} reports empirical
sizes of the proposed test and of the CZZ test for the preceding three
distributions.
We observe from Table~\ref{table2} that the sizes of both tests are roughly the
same, when the underlying variables are normally or Bernoulli
distributed. It seems that the CZZ test looks better for skewed data,
for example, gamma distribution. We believe additional corrections such as the
Edgeworth expansion will be helpful, which is beyond the scope of this
paper. 
However, our test still performs well for skewed data if $p\gg n$.

Table~\ref{table3} to Table~\ref{table5} summarize the empirical powers of the proposed tests
as well as those of the CZZ test for both the diagonal and the banded
covariance matrix. Table~\ref{table3} assumes the underlying variables are
normally distributed while Tables~\ref{table4} and \ref{table5} assume the central gamma and
the central bernoulli random variables, respectively. For the diagonal
covariance matrix, we observe that the proposed test consistently
outperforms the CZZ test for all types of distributions, especially for
``small'' $n$. For example, when $n=20$, even $n=40,60,80$ for $\nu
=0.08$, the CZZ test results in power ranging from 0.2--0.8, while our
test still gains very satisfying power exceeding 0.932.

For the banded covariance matrix, we observe an interesting phenomenon.
Our test seems to be more sensitive to the dimension $p$. When
$p=600,1500,3000$, the power of our test is not that good for small
$v_2$ ($=0.4$). Fortunately, when $p=5500,8000,10\,000$, the performance is
much better, where the power is one or close to one. Similar results
are also observed for $v_2=0.8$. We also note that large $v_2$
outperforms smaller $v_2$ because when $v_2$ becomes larger, the
corresponding covariance matrix becomes more ``different'' from the
identity matrix. As for the CZZ test, its power is mainly affected by
$n$. But generally speaking, our test gains better power than the CZZ
test for extremely larger $p$ and small $n$.


\subsection{Empirical studies}\label{empirical}

As empirical applications, we consider two classic datasets: the colon
data of Alon \textit{et al.} \cite{alon99} and the leukemia data of
Golub \textit{et al.} \cite{Golub99}. Both datasets are publicly available on the web
site of Tatsuya Kubokawa: \url{http://www.tatsuya.e.u-tokyo.ac.jp/}. Such
data were used in Fisher \cite{Fisher12} as well. The
sample sizes and
dimensions $(n,p)$ of the colon data and the leukemia data are
$(62,2000)$ and $(72,3571)$, respectively. Simulations show that these
two datasets have zero mean ($10^{-8}$ to $10^{-11}$) and unit
variance. Therefore, we consider the hypothesis test in \eqref{k1} by
using the test statistic $L_n$ in \eqref{o4}. The computed values are
$L_n=33\,933.7$ for the colon data and $L_n=60\,956$ for the leukemia data.
It is also interesting to note that the statistic values of Fisher \cite{Fisher12} are $6062.642$ for the colon data and
$6955.651$ for the
leukemia data when testing the identity hypothesis. Also, the
statistics of Fisher \cite{Fisher12} and $L_n$ in
\eqref{o4} are both
asymptotic normality (standard normal). As in Fisher
\cite{Fisher12}, we
conclude that $p$-values of the test statistics are zero which shows
evidence to reject the null hypothesis. This is consistent with
Fisher's \cite{Fisher12} conclusion for these two datasets.


\section{Truncation and strategy for the proof of Theorem~\texorpdfstring{\protect\ref{thm2}}{1.1}}\label{TandS}

In the rest of the paper, we use $K$ to denote a constant which may
take different values at different places. The notation $\mathrm{o}_{L_p}(1)$
stands for a term converging to zero in $L_p$ norm; $\xrightarrow
{\mathrm{a.s.}}$ means ``convergence almost surely to''; $\xrightarrow{\mathrm{i.p.}}$
means ``convergence in probability to''.

\subsection{Truncation}\label{truncation}

In this section, we truncate the underlying random variables as in
Pan and Gao \cite{PG12}. Choose $\delta_n$ satisfying
%
\begin{equation}
\label{c1} \lim_{n\to\infty}\delta_n^{-4}E|X_{11}|^4I
\bigl(|X_{11}|>\delta_n\sqrt [4]{np}\bigr)=0, \qquad
\delta_n\downarrow0, \delta_n\sqrt[4]{np}\uparrow
\infty.
\end{equation}
In what follows, we will use $\delta$ to represent $\delta_n$ for
convenience. We first truncate the variables $\hat
{X}_{ij}=X_{ij}I(|X_{ij}|<\delta\sqrt[4]{np})$ and then normalize it
as $\tilde{X}_{ij}=(\hat{X}_{ij}-E\hat{X}_{ij})/\sigma$, where
$\sigma$ is the standard deviation of $
\hat{X}_{ij}$. Let $\hat\bbX=(\hat{X}_{ij})$ and $\tilde\bbX
=(\tilde{X}_{ij})$. Define $\hat{\bbA},\tilde{\bbA}$ and $\hat
{G}_n(f),\tilde{G}_n(f)$ similarly by means of (\ref{v1}) and (\ref
{z22}), respectively. We then have
\[
P(\bbA\neq\hat{\bbA})\leq npP\bigl(|X_{11}|\geq\delta\sqrt [4]{np}
\bigr)\leq K\delta^{-4}E|X_{11}|^4I
\bigl(|X_{11}|>\delta\sqrt[4]{np}\bigr)=\mathrm{o}(1).
\]

It follows from (\ref{c1}) that
\begin{eqnarray*}
\bigl|1-\sigma^2\bigr|&\leq&2\bigl|EX_{11}^2I
\bigl(|X_{11}>\delta\sqrt[4]{np}|\bigr)\bigr| \\
&\leq&2(np)^{-1/2}
\delta^{-2}E|X_{11}|^4I\bigl(|X_{11}|>
\delta\sqrt [4]{np}\bigr)=\mathrm{o} \bigl((np)^{-1/2} \bigr)
\end{eqnarray*}
and
\[
|E\hat{X}_{11}|\leq\delta ^{-3}(np)^{-3/4}E|X_{11}|^4I
\bigl(|X_{11}|>\delta\sqrt[4]{np}\bigr)=\mathrm{o} \bigl((np)^{-3/4}
\bigr).
\]
Therefore
\begin{eqnarray*}
&&E\tr(\tilde{\bbX}-\hat{\bbX})^{T}(\tilde{\bbX}-\hat{\bbX})\leq \sum
_{i,j}E|\hat{X}_{ij}-\tilde{X}_{ij}|^2\\
&&\quad
\leq Kpn \biggl(\frac{(1-\sigma)^2}{\sigma^2}E|\hat{X}_{11}|^2+
\frac
{1}{\sigma^2}|E\hat{X}_{ij}|^2 \biggr) =\mathrm{o}(1)
\end{eqnarray*}
and
\[
E\tr\hat{\bbX}^T\hat{\bbX}\leq\sum_{i,j}E|
\hat{X}_{ij}|^2\leq Knp,\qquad  E\tr\tilde{\bbX}\tilde{
\bbX}^T\leq\sum_{i,j}E|\tilde
{X}_{ij}|^2\leq Knp.
\]

Recalling that the notation $\lambda_j(\cdot)$ represents the $j$th
largest eigenvalue, we then have $\lambda_j(\bbX^T\bbX)=\sqrt {np}\lambda_j(\bbA)+p$. Similar equalities also hold if $\bbX,\bbA$
are replaced by $\hat{\bbX},\hat{\bbA}$ or $\tilde{\bbX},\tilde
{\bbA}$. Consequently, applying the argument used in Theorem~11.36 in
Bai and Silverstein \cite{BSbook} and Cauchy--Schwarz's inequality, we have
\begin{eqnarray*}
E\bigl\llvert \tilde{G}_n(f)-\hat{G}_n(f)
\bigr\rrvert &\leq&\sum_{j=1}^{n}E\bigl|f\bigl(
\lambda_j(\hat{\bbA})\bigr)-f\bigl(\lambda_j(\tilde{
\bbA})\bigr)\bigr|
\\
&\leq& K_f \sum_{j=1}^{n}E\bigl|
\lambda_j(\hat{\bbA})-\lambda_j(\tilde {\bbA})\bigr|\leq
\frac{K_f}{\sqrt{np}}\sum_{j=1}^{n}E\bigl|
\lambda_j\bigl(\hat {\bbX}^T\hat{\bbX}\bigr)-
\lambda_j\bigl(\tilde{\bbX}^T\tilde{\bbX}\bigr)\bigr|
\\
&\leq&\frac{K_f}{\sqrt{np}}E \bigl[\tr(\tilde{\bbX}-\hat{\bbX })^T(
\tilde{\bbX}-\hat{\bbX})\cdot2 \bigl(\tr\hat{\bbX}^T\hat {\bbX}+\tr
\tilde{\bbX}^T\tilde{\bbX} \bigr) \bigr]^{1/2}
\\
&\leq&\frac{2K_f}{\sqrt{np}} \bigl[E\tr(\tilde{\bbX}-\hat{\bbX })^{T}(
\tilde{\bbX}-\hat{\bbX})\cdot \bigl(E\tr\hat{\bbX}^T\hat {\bbX} + E\tr
\tilde{\bbX}^T\tilde{\bbX} \bigr) \bigr]^{\sfrac{1}{2}}
\\
&=&\mathrm{o}(1), %
\end{eqnarray*}
where $K_f$ is a bound on $|f'(x)|$. Thus, the weak convergence of
$G_n(f)$ is not affected if we replace the original variables $X_{ij}$
by the truncated and normalized variables $\tilde{X}_{ij}$. For
convenience, we still use $X_{ij}$ to denote $\tilde{X}_{ij}$, which
satisfies the following additional assumption (c):
\begin{enumerate}[(c)]
\item[(c)] The underlying variables satisfy
\[
|X_{ij}|\leq\delta\sqrt[4]{np},\qquad  EX_{ij}=0,\qquad
EX_{ij}^2=1, \qquad EX_{ij}^4=
\nu_4+\mathrm{o}(1),
\]
where $\delta=\delta_n$ satisfies $\lim_{n\to\infty}\delta
_n^{-4}E|X_{11}|^4I(|X_{11}|>\delta_n\sqrt[4]{np})=0$, $\delta
_n\downarrow0$, and\linebreak[4]  $\delta_n\sqrt[4]{np}\uparrow\infty$.
\end{enumerate}

For any $\epsilon>0$, define the event $F_n(\epsilon)=\{\max_{j\leq
n}|\lambda_j(\bbA)|\geq2+\epsilon\}$ where $\bbA$ is defined by
the truncated and normalized variables satisfying assumption (c). By Theorem~2 in Chen and Pan \cite{CP12}, for any $\ell>0$
%
\begin{equation}
\label{qstar} P\bigl(F_n(\epsilon)\bigr)=\mathrm{o}
\bigl(n^{-\ell}\bigr).
\end{equation}
Here we would point out that the result regarding the minimum
eigenvalue of $\bbA$ can be obtained similarly by investigating the
maximum eigenvalue of $-\bbA$.

\subsection{Strategy of the proof}\label{strategy}
We shall follow the strategy of Bai and Yao \cite{BY05}.
Specifically speaking,
assume that $u_0,v$ are fixed and sufficiently small so that $\varsigma
\subset\mathscr{S}$ (see the definition in the introduction), where
$\varsigma$ is the contour formed by the boundary of the rectangle
with $(\pm u_0, \pm \mathrm{i}v)$ where $u_0>2,0<v\leq1$. By Cauchy's integral
formula, with probability one,
\[
G_n(f)=-\frac{1}{2\pi \mathrm{i}}\oint_{\varsigma}f(z)n
\bigl[m_n(z)-m(z)-\mathcal{X}_n\bigl(m(z)\bigr) \bigr]
\dd z,
\]
where $m_n(z),m(z)$ denote the Stieltjes transform of $F^{\bbA}(x)$
and $F(x)$, respectively.

Let
\[
M_n(z)=n \bigl[m_n(z)-m(z)-\mathcal{X}_n
\bigl(m(z)\bigr) \bigr],\qquad  z\in\varsigma.
\]
For $z\in\varsigma$, write $M_n(z)=M_n^{(1)}(z)+M_n^{(2)}(z)$ where
\[
M_n^{(1)}(z)=n\bigl[m_n(z)-Em_n(z)
\bigr],\qquad  M_n^{(2)}(z)=n \bigl[Em_n(z)-m(z)-
\mathcal{X}_n\bigl(m(z)\bigr) \bigr].
\]

Split the contour $\varsigma$ as the union of $\varsigma_u, \varsigma
_l, \varsigma_r,\varsigma_0$ where $\varsigma_l=\{z=-u_0+\mathrm{i}v,\xi
_nn^{-1}<|v|<v_1\},\varsigma_r=\{z=u_0+\mathrm{i}v,\xi_nn^{-1}<|v|<v_1\}
,\varsigma_0=\{z=\pm u_0+\mathrm{i}v,|v|\leq\xi_nn^{-1}\}$ and $\varsigma
_u=\{z=u\pm \mathrm{i}v_1,|u|\leq u_0\}$ and where $\xi_n$ is a slowly varying
sequence of positive constants and $v_1$ is a positive constant which
is independent of $n$. Throughout this paper, let $\mathbb{C}_1=\{
z\dvtx z=u+\mathrm{i}v, u\in[-u_0,u_0], |v|\geq v_1\}$.

\begin{proposition}\label{pro1}
Under assumptions \textup{(b1)}, \textup{(c)}, the empirical process
$\{M_n(z),z\in\mathbb{C}_1\}$ converges weakly to a Gaussian process
$\{M(z),z\in\mathbb{C}_1\}$ with the mean function
\begin{equation}
\label{x3} \Delta(z)=0
\end{equation}
and the covariance function
%
\begin{equation}
\label{x7} \Lambda(z_1,z_2)=m'(z_1)m'(z_2)
\bigl[\nu_4-3+2 \bigl(1-m(z_1)m(z_2)
\bigr)^{-2} \bigr].
\end{equation}
\end{proposition}

%

As in Bai and Yao \cite{BY05}, the process of $\{M(z),z\in
\mathbb{C}_1\}$ can be
extended to $\{M(z),\Re(z)\notin[-2,2]\}$ due to the facts that (i)
$M(z)$ is symmetric, for example, $M(\bar{z})=\overline{M(z)}$; (ii) the mean
and the covariance function of $M(z)$ are independent of $v_1$ and they
are continuous except for $\Re(z)\notin[-2,2]$. By Proposition~\ref
{pro1} and the continuous mapping theorem,
%
\begin{eqnarray*}
-\frac{1}{2\pi \mathrm{i}}\int_{\varsigma_u}f(z)M_n(z)\dd z
\xrightarroww {d}-\frac{1}{2\pi \mathrm{i}}\int_{\varsigma_u}f(z)M(z)\dd z.
\end{eqnarray*}
Thus, to prove Theorem~\ref{thm2}, it is also necessary to prove the
following proposition.

\begin{proposition}\label{pro3}
Let $z\in\mathbb{C}_1$. Under assumptions \textup{(b1)},
\textup{(c)}, there exists some event $U_n$ with $P(U_n)\to0$, as $n\to
\infty$, such that
%
\begin{eqnarray}
\label{y11} \lim_{v_1\downarrow0}\limsup_{n\to\infty}E
\biggl\llvert \int_{\bigcup
_{i=l,r,0}\varsigma_i}M_n^{(1)}(z)I
\bigl(U_n^c\bigr)\dd z\biggr\rrvert ^2&=&0,
\\
\label{y1} \lim_{v_1\downarrow0}\limsup_{n\to\infty}\biggl
\llvert \int_{\bigcup
_{i=l,r,0}\varsigma_i}EM_n(z)I\bigl(U_n^c
\bigr)\dd z\biggr\rrvert &=&0
\end{eqnarray}
%
and
%
\begin{equation}
\label{y2} \lim_{v_1\downarrow0}E\biggl\llvert \int
_{\varsigma_i}M^{(1)}(z)\dd z\biggr\rrvert ^2=0,\qquad
\lim_{v_1\downarrow0}E\biggl\llvert \int_{\varsigma_i}M(z)\dd
z\biggr\rrvert ^2=0.
\end{equation}
\end{proposition}

Since $E|M^{(1)}(z)|^2=\Lambda(z,\bar{z})$ and $E|M(z)|^2=\Lambda
(z,\bar{z})+|EM(z)|^2$, (\ref{y2}) can be easily obtained from
Proposition~\ref{pro1}.
For $i=0$, if we choose $U_n= F_n(\epsilon)$ with the $\epsilon
=(u_0-2)/2$, then when $U_n^c$ happens, $\forall z\in\varsigma_0$, we
have $|m_n(z)|\leq2/(u_0-2)$ and $|m(z)|\leq1/(u_0-2)$. Thus
\[
\biggl|\int_{\varsigma_0}M_n^{(1)}(z)I
\bigl(U_n^c\bigr)\dd z \biggr|\leq n \biggl(\frac
{4}{u_0-2}
\biggr)^2\|\varsigma_0\|\leq\frac{4\xi_n}{(u_0-2)^2},
\]
where $\|\varsigma_0\|$ represents the length of $\varsigma_0$. Furthermore,
\[
\biggl|\int_{\varsigma_0}M_n(z)I\bigl(U_n^c
\bigr)\dd z \biggr|\leq n \biggl(\frac
{2}{u_0-2}+\frac{1}{u_0-2}+K
\frac{n}{p} \biggr)^2\|\varsigma_0\|.
\]
These imply that (\ref{y1}) and (\ref{y11}) are true for $z\in
\varsigma_0$ by noting that $\xi_n\rightarrow0$ as $p\rightarrow
\infty$.

Sections~\ref{M1} and \ref{M2} are devoted to the proof of
Proposition~\ref{pro1}. The main steps are summarized in the following:
\begin{itemize}
\item According to Theorem~8.1 in Billingsley \cite{bili68},
to establish the
convergence of the process $\{M_n(z),z\in\mathbb{C}_1\}$, it suffices
to prove the finite-dimensional convergence of the random part
$M_n^{(1)}(z)$ and its tightness, and the convergence of the non-random
part $M_n^{(2)}(z)$.
\item For the random part $M_n^{(1)}(z)$, we rewrite it in terms of a
martingale expression 
so that we may apply the central limit theorem of martingales to find
its asymptotic mean and covariance.
%
\item For the non-random part $M_n^{(2)}(z)$, by the formula of the
inverse of a matrix and the equation satisfied by $m(z)$ we develop an
equation for $(Em_n(z)-m(z))$. Based on it, we then find its limit
under assumptions $n/p\to0$ and $n^3/p=\mathrm{O}(1)$ for Theorem~\ref{thm2}
and Corollary~\ref{cor1}, respectively.
\end{itemize}

Section~\ref{proofPro3} uses Lemma~\ref{lem7} below to finish the
proofs of \eqref{y11} for $i=l,r$ so that the proof of Proposition~\ref{pro3} is completed. Section~\ref{calculation} uses
Bai and Yao's \cite{BY05} asymptotic mean and covariance function to
conclude the proof
of Theorem~\ref{thm2}.

\section{Preliminary results}\label{lemmas}

This section is to provide simplification of $M_n^{(1)}(z)$ and some
useful lemmas needed to prove Proposition~\ref{pro1}.

\subsection{Simplification of $M_n^{(1)}(z)$}
\label{simplification}

The aim of this subsection is to simplify $M_n^{(1)}(z)$ so that
$M_n^{(1)}(z)$ can be written in the form of martingales. Some moment
bounds are also proved.

Define $\bbD=\bbA-z\bbI_n$. Let $\mathbf{s}_k$ be the $k$th column of
$\bbX$ and $\bbX_k$ be a $p\times(n-1)$ matrix constructed from
$\bbX$ by deleting the $k$th column. We then similarly define $\bbA
_k=\frac{1}{\sqrt{np}}(\bbX_k^T\bbX_k-p\bbI_{n-1})$ and $\bbD
_k=\bbA_k-z\bbI_{n-1}$. The $k$th diagonal element of $\bbD$ is
$a_{kk}^{\diagg}=\frac{1}{\sqrt{np}}(\mathbf{s}_k^T\mathbf
{s}_k-p)-z$ and the
$k$th row of $\bbD$ with the $k$th element deleted is $\bbq_k^T=\frac
{1}{\sqrt{np}}\mathbf{s}_k^T\bbX_k$. The Stieltjes transform of
$F^{\bbA
}$ has the form $m_n(z)=\frac{1}{n}\tr\bbD^{-1}$. The limiting
Stieltjes transform $m(z)$ satisfies
%
\begin{equation}
\label{q4} m(z)=-\frac{1}{z+m(z)}, \qquad \bigl|m(z)\bigr|\leq1
\end{equation}
(one may see Bai and Yao \cite{BY05}).

Define the $\sigma$-field $\mathcal{F}_k=\sigma(\mathbf
{s}_1,\mathbf{s}_2,\ldots,\mathbf{s}_{k})$ and the conditional
expectation $E_k(\cdot
)=E(\cdot|\mathcal{F}_k)$. By the matrix inversion formula, we have
(see (3.9) of Bai \cite{Bai93})
%
\begin{equation}
\label{h11} \tr\bigl(\bbD^{-1}-\bbD_k^{-1}
\bigr)=-\frac{(1+\bbq_k^{T} \bbD_k^{-2}\bbq
_k)}{-a_{kk}^{\diagg}+\bbq_k^T\bbD_k^{-1}\bbq_k}.
\end{equation}
We then obtain
%
\begin{eqnarray}
\label{u4} M_n^{(1)}(z)&=&\tr\bbD^{-1}-E\tr
\bbD^{-1}=\sum_{k=1}^{n}(E_{k}-E_{k-1})
\tr\bigl(\bbD^{-1}-\bbD_k^{-1}\bigr) = \sum
_{k=1}^{n}\varrho_k
\\
&=&(E_{k}-E_{k-1})\iota_k-E_{k}
\kappa_k,
\end{eqnarray}
where
\begin{eqnarray*}
\varrho_k&=&-(E_{k}-E_{k-1})
\beta_k \bigl(1+\bbq_k^T
\bbD_k^{-2}\bbq _k \bigr),
\\
\iota_k&=&-\beta_k^{\trr}\beta_k
\eta_k \bigl(1+\bbq_k^T\bbD
_k^{-2}\bbq_k \bigr),
\\
\eta_k&=&\frac{1}{\sqrt{np}}\bigl(\mathbf{s}_k^T
\mathbf{s}_k-p\bigr)-\gamma_{k1},\qquad  \beta_k=
\frac{1}{-a_{kk}^{\diagg}+\bbq_k^T\bbD_k^{-1}\bbq_k},
\\
\beta_k^{\trr}&=&\frac{1}{z+\sklfrac{1}{(np)}\tr\bbM_k^{(1)}},\qquad  \bbM
_k^{(s)}=\bbX_k\bbD_k^{-s}
\bbX_k^T, s=1,2,
\\
\gamma_{ks}& =& \bbq_k^T\bbD_k^{-s}
\bbq_k-(np)^{-1}\tr\bbM_k^{(s)},\qquad
\kappa_k=\beta_k^{\trr}\gamma_{k2}.
\end{eqnarray*}
In the above equality, $\varrho_k$ is obtained by (\ref{h11}) and the
last equality uses the facts that
%
\begin{equation}
\label{h1} \beta_k=\beta^{\trr}_k+
\beta_k\beta^{\trr}_k\eta_k
\end{equation}
and
\[
(E_{k}-E_{k-1}) \biggl[\beta_k^{\trr}
\biggl(1+ \frac{1}{np}\tr\bbM _k^{(2)} \biggr)
\biggr]=0,\qquad  E_{k-1}\kappa_k=0.
\]
We remind the readers that the variable $z$ has been dropped from the
expressions such as $\bbD^{-1},\bbD_k^{-1}, \beta_k$, $\gamma_{ks}$
and so on. When necessary, we will also indicate them as $\bbD
^{-1}(z),\bbD_k^{-1}(z),\allowbreak  \beta_k(z)$, $\gamma_{ks}(z)$, etc.

We next provide some useful bounds. It follows from the definitions of
$\bbD$ and $\bbD_k$ that
%
\begin{eqnarray}
\label{a0} %
\bbD^{-1}\bbX^T\bbX&=&p
\bbD^{-1}+\sqrt{np}\bigl(\bbI_n+z\bbD^{-1}\bigr),
\nonumber\\[-8pt]\\[-8pt]
\bbD_k^{-1}\bbX_k^T
\bbX_k&=&p\bbD_k^{-1}+\sqrt{np}\bigl(\bbI
_{n-1}+z\bbD_k^{-1}\bigr).\nonumber %
\end{eqnarray}
Since the eigenvalues of $\bbD^{-1}$ have the form $1/(\lambda_j(\bbA
)-z)$, $\|\bbD^{-1}\|\leq1/v_1$ and similarly $\|\bbD_k^{-1}\|\leq1/v_1$.
From Theorem~11.4 in Bai and Silverstein \cite{BSbook}, we note that
$-\beta_k(z)$ is the
$k$th diagonal element of $\bbD^{-1}$ so that $|\beta_k|\leq1/v_1$.
Moreover, considering the imaginary parts of $1/\beta_k^{\trr}$ and
$1/\beta_k$ and by (\ref{a0}) we have
%
\begin{equation}
\label{h2} \bigl|\beta_k^{\trr}\bigr|\leq1/v_1,\qquad  \biggl|1+
\frac{1}{np}\tr\bbM_k^{(s)}\biggr| \leq
\bigl(1+1/v_1^{2s}\bigr),\qquad s=1,2
\end{equation}
and
%
\begin{equation}
\label{h3} \bigl\llvert \bigl(1+\bbq_k^T
\bbD_k^{-2}\bbq_k \bigr)\beta_k
\bigr\rrvert \leq\frac{1+\bbq_k^T\bbD_k^{-1}\overline{\bbD}_k^{-1}\bbq
_k}{v_1(1+\bbq_k^T\bbD_k^{-1}\overline{\bbD}_k^{-1}\bbq_k)}= 1/v_1.
\end{equation}

Applying (\ref{h1}), we split $\iota_k$ as
\begin{eqnarray*}
\iota_k &=&- \biggl(1+\frac{1}{np}\tr
\bbM_k^{(2)} \biggr) \bigl(\beta ^{\trr}_k
\bigr)^2\eta_k-\gamma_{k1}\bigl(
\beta^{\trr}_k\bigr)^2\eta_k-
\biggl(1+\frac
{1}{np}\bbq_k^T\bbD_k^{-2}
\bbq_k \biggr) \bigl(\beta^{\trr}_k
\bigr)^2\beta _k\eta_k^2
\\
& =& \iota_{k1}+\iota_{k2}+\iota_{k3}.
\end{eqnarray*}
As will be seen, $\iota_{k1},\iota_{k2}$ could be negligible by Lemma~\ref{lem1} below.


By Lemma~\ref{lem1}, \eqref{h2} and \eqref{h3}, we have
\begin{eqnarray*}
E \Biggl|\sum_{k=1}^{n}(E_{k}-E_{k-1})
\iota_{k3} \Biggr|^2\leq\sum_{k=1}^{n}E
\biggl\llvert \biggl(1+\frac{1}{np}\mathbf{s}_k^T
\bbM _k^{(2)}\mathbf{s}_k \biggr) \bigl(
\beta_k^{\trr}\bigr)^2\beta_k
\eta_k^2\biggr\rrvert ^2\leq K
\delta^4,
\end{eqnarray*}
and that
\begin{eqnarray*}
E \Biggl|\sum_{k=1}^{n}(E_{k}-E_{k-1})
\iota_{k2} \Biggr|^2 \leq\sum_{k=1}^{n}E
\bigl\llvert \gamma_{k1}\bigl(\beta_k^{\trr}
\bigr)^2\eta_k\bigr\rrvert ^2\leq K\sum
_{k=1}^{n} \bigl(E|\gamma_{k1}|^4E|
\eta_k|^4 \bigr)^{1/2} \leq \frac{Kn}{p}+K
\delta^2.
\end{eqnarray*}
Therefore, $M_n^{(1)}(z)$ is simplified as
%
\begin{eqnarray}
\label{a3} %
M_n^{(1)}(z)&=&\sum
_{k=1}^{n}E_{k} \biggl[- \biggl(1+
\frac{1}{np}\tr \bbM _k^{(2)} \biggr) \bigl(
\beta_k^{\trr}\bigr)^2 \eta_k-
\kappa_k \biggr]+\mathrm{o}_{L_2}(1)\nonumber
\\[-8pt]\\[-8pt]
& =& \sum_{k=1}^{n}E_{k}\bigl(
\alpha_k(z)\bigr)+\mathrm{o}_{L_2}(1),\nonumber %
\end{eqnarray}
where $\alpha_k(z)$ represents the term in the square bracket. Thus,
to prove finite-dimensional convergence of $ M_n^{(1)}(z),z\in\mathbb
{C}_1$ we need only consider the sum
%
\begin{equation}
\label{h12} \sum_{j=1}^l
a_j \sum_{k=1}^{n}E_{k}
\bigl(\alpha_k(z_j)\bigr)= \sum
_{k=1}^{n}\sum_{j=1}^l
a_jE_{k}\bigl(\alpha_k(z_j)
\bigr),
\end{equation}
where $a_1,\ldots,a_l$ are complex numbers and $l$ is any positive integer.

\subsection{Useful lemmas}
The aim of this subsection is to provide some useful lemmas.

%
\begin{lem}\label{lem1}
Let $z\in\mathbb{C}_1$. Under assumptions \textup{(b1)},
\textup{(c)}, we have
%
\begin{eqnarray}
\label{h10} E| \gamma_{ks}|^2&\leq& Kn^{-1},\qquad  E|
\gamma_{ks}|^4\leq K \biggl(\frac
{1}{n^2}+
\frac{n}{p^2}+\frac{1}{np} \biggr),
\\
\label{x9} E|\eta_k|^2&\leq& Kn^{-1},\qquad  E|
\eta_k|^4\leq K\frac{\delta^4}{n}+K \biggl(
\frac{1}{n^2}+\frac{p}{n^2}+\frac{1}{np} \biggr).
\end{eqnarray}
\end{lem}

\begin{pf} From Lemma~5 in Pan and Zhou \cite{PZ11}, we obtain
%
\begin{equation}
\label{h9} E\bigl|\mathbf{s}_k^T\bbH\mathbf{s}_k-
\tr\bbH\bigr|^4\leq K \bigl(EX_{11}^4
\bigr)^2E(\tr\bbH\bbH)^2\leq KE \bigl(\tr
\bbM_k^{(s)}\overline{\bbM }_k^{(s)}
\bigr)^2\leq Kn^2p^4,
\end{equation}
where $\bbH=\bbM_k^{(s)}-\diag(a_{11}^{(s)},\ldots,a_{nn}^{(s)})$ and
$a_{jj}^{(s)}$ is the $j$th diagonal element of the matrix $\bbM
_k^{(s)}$. To get the third inequality in (\ref{h9}), by (\ref{a0})
and the uniform bound for $\|\bbD_k^{-1}\|$, we obtain
%
\begin{eqnarray}
\label{i10} %
\bigl|\tr\bbM_k^{(s)}\overline{
\bbM}_k^{(s)}\bigr|&=&\bigl|\tr\bbD_k^{-s}\bbX
_k^T\bbX_k\overline{\bbD}_k^{-s}
\bbX_k^T\bbX_k\bigr|\leq\frac
{n}{v_1^{2(s-1)}}\bigl\|
\bbD_k^{-1}\bbX_k^T
\bbX_k\bigr\|^2
\nonumber\\[-8pt]\\[-8pt]
&\leq&\frac{n}{v_1^{2(s-1)}}\bigl\|p\bbD_k^{-1}+\sqrt{np}\bigl(\bbI
_{n-1}+z\bbD_k^{-1}\bigr)\bigr\|^2\leq
\frac{Kn^2p^4}{v_1^{2s}}.\nonumber %
\end{eqnarray}
Let $\mathbb{E}_j(\cdot)=E(\cdot|X_{1k},X_{2k},\dots
,X_{jk}),j=1,\dots,p$. Since $\{X_{jk}\}_{j=1}^{k}$ are independent of
$a_{jj}^{(s)}$, $(X_{jk}^2-1)a_{jj}^{(s)}=(\mathbb{E}_j-\mathbb
{E}_{j-1})(X_{jk}^2-1)a_{jj}^{(s)}$. By Burkholder's inequality and
assumption (c)
%
\begin{eqnarray}
\label{h8} %
E\Biggl|\sum_{j=1}^p
\bigl(X_{jk}^2-1\bigr)a_{jj}^{(s)}\Biggr|^4&=&E
\Biggl|\sum_{j=1}^p(\mathbb{E}_j-
\mathbb{E}_{j-1}) \bigl(X_{jk}^2-1
\bigr)a_{jj}^{(s)} \Biggr|^4
\nonumber\\
&\leq&KE \Biggl(\sum_{j=1}^nE|X_{11}|^4\bigl|a_{jj}^{(s)}\bigr|^2
\Biggr)^2+K\sum_{j=1}^pE|X_{11}|^8E|a_{jj}|^4\\
&\leq& Kn^5p^2+n^3p^3,\nonumber %
\end{eqnarray}
where we use the fact that, with $\mathbf{w}_j^T$ being the $j$th row of
$\bbX_k$,
%
\begin{eqnarray}\label{h7}
E\bigl|a_{jj}^{(s)}\bigr|^4&=&E\bigl|\ddot{\bbe}_j^T
\bbX_k\bbD_k^{-s}\bbX_k^T
\ddot {\bbe}_j\bigr|^4\nonumber\\[-8pt]\\[-8pt]
&=&E\bigl|\mathbf{w}_j^T
\bbD_k^{-s}\mathbf{w}_j\bigr|^4\leq
v_1^{-4s}E\bigl\|\mathbf{w}_j^T
\bigr\|^8\leq Kn^4+Kn^2p.\nonumber
\end{eqnarray}
Here for $j=1,\dots,p$, $\ddot{\bbe}_j$ denotes the $p$-dimensional
unit vector with the $j$th element being $1$ and all the remaining
being zero. It follows from (\ref{h9}) and (\ref{h8}) that
\begin{eqnarray*}
E| \gamma_{ks}|^4&\leq& \frac{K}{n^4p^4} E\Biggl|\sum
_{j=1}^p\bigl(X_{jk}^2-1
\bigr)a_{jj}^{(s)}\Biggr|^4+\frac{K}{n^4p^4}E\bigl|\mathbf
{s}_k^T\bbH \mathbf{s}_k-\tr
\bbH\bigr|^4
\nonumber
\\
&\leq&K \biggl(\frac{1}{n^2}+\frac{n}{p^2}+\frac{1}{np} \biggr).
\end{eqnarray*}
Moreover, applying Lemma~8.10 in Bai and Silverstein \cite{BSbook},
we have
\begin{eqnarray*}
E|\eta_k|^4\leq\frac{K}{n^2p^2}E\bigl|
\mathbf{s}_k^T\mathbf{s}_k-n\bigr|^4+KE\bigl|
\gamma_{k1}(z)\bigr|^4\leq K\frac{\delta^4}{p}+K \biggl(
\frac
{1}{n^2}+\frac{p}{n^2}+\frac{1}{np} \biggr). %
\end{eqnarray*}
The bounds of the absolute second moments for $\gamma_{ks},\eta_k$
follow from a direct application of Lemma~8.10 in Bai and Silverstein
\cite{BSbook}, (\ref
{a0}) and the uniform bound for $\|\bbD_k^{-1}\|$.
\end{pf}


When $z\in\varsigma_l\cup\varsigma_r$, the spectral norm of $\bbD
^{-1}(z)$ as well as the quantities in \eqref{h2} or Lemma~\ref{lem1}
are unbounded. In order to prove Lemma~\ref{lem4}, we will establish
the bounds similar to those in \eqref{h2} and in Lemma~\ref{lem1} for
$z\in\varsigma_l\cup\varsigma_r$ below.

Let the event $U_n=\{\max_{j\leq n}|\lambda_j(\bbA)|\geq u_0/2+1\}$
and $U_{nk}=\{\max_{j\leq n}|\lambda_j(\bbA_k)|\geq1+u_0/2\}$.
The Cauchy interlacing theorem ensures that
%
\begin{equation}
\label{s6} \lambda_1(\bbA)\geq\lambda_1(
\bbA_k)\geq\lambda_2(\bbA)\geq \lambda_2(
\bbA_k)\geq\cdots\geq\lambda_{n-1}(\bbA_k)\geq
\lambda _n(\bbA).
\end{equation}
Thus, $U_{nk}\subset U_n$. By (\ref{qstar}) for any $\ell>0$
%
\begin{equation}
\label{s1} P(U_{nk})\leq P(U_n)=\mathrm{o}
\bigl(n^{-\ell}\bigr).
\end{equation}
We claim that
%
\begin{eqnarray}
\label{s0} \hspace*{-30pt}\max\bigl\{\bigl\|\bbD^{-1}(z)\bigr\|,\bigl\|\bbD^{-1}(z)\bigr\|,|
\beta_k|\bigr\}&\leq&\xi_n^{-1}n;
\\
\label{s7} \frac{I(U_n^c)}{|\lambda_j(\bbA)-z|}&\leq& K,\qquad  j=1,2,\ldots,n, \nonumber\\[-8pt]\\[-8pt]
 \frac{
I(U_{nk}^c)}{|\lambda_{j}(\bbA_k)-z|}&\leq&
K, \qquad i=1,2,\ldots,(n-1);\nonumber
\\
\label{add8} \bigl\|\bbD^{-1}(z)\bigr\|I\bigl(U_n^c
\bigr)&\leq&2/(u_0-2),\qquad  \bigl\|\bbD_k^{-1}(z)\bigr\| I
\bigl(U_{nk}^c\bigr)\leq2/(u_0-2).
\end{eqnarray}
Indeed, the quantities in (\ref{s0}) are bounded due to $|1/\Im
(z)|\leq\xi_n^{-1}n$ while (\ref{s7}) holds because
$I(U_n^c)/|\lambda_j(\bbA)-z|$ (or $I(U_n^c)/|\lambda_j(\bbA
_k)-z|$) is bounded by $v_1^{-1}$ when $z\in\varsigma_u$ and bounded
by $2/(u_0-2)$ when $z\in\varsigma_l \cup\varsigma_r$. The
estimates in (\ref{add8}) hold because of the eigenvalues of $\bbD
^{-1}I(U_n^c)$ (or $\bbD^{-1}I(U_n^c)$) having the form
$I(U_n^c)/(\lambda_j(\bbA)-z)$ (or $I(U_n^c)/(\lambda_j(\bbA_k)-z)$).

%
\begin{lem}\label{lem5}
Let $z\in\varsigma_n$. The following bound
%
\begin{equation}
\label{add1} |\beta_k|I\bigl(U_n^c\bigr)
\leq K,
\end{equation}
holds.
\end{lem}

\begin{pf} In view of (\ref{h11}), to prove (\ref{add1}), we
need to find an upper bound for $ |\tr\bbD^{-1}-\tr\bbD
_k^{-1}
|I(U_n^c)$ and a lower bound for $|1+\bbq_k^T\bbD_k^{-2}\bbq
_k|I(U_n^c)$. It follows from (\ref{s7}) and (\ref{s6}) that
%
\begin{eqnarray}
\label{s8} %
\bigl\llvert \tr\bbD^{-1}-\tr
\bbD_k^{-1}\bigr\rrvert I\bigl(U_n^c
\bigr)&\leq&\Biggl\llvert \sum_{j=1}^{n}
\frac{1}{\lambda_j(\bbA)-z}-\sum_{j=1}^{n-1}
\frac
{1}{\lambda_{j}(\bbA_k)-z}\Biggr\rrvert I\bigl(U_n^c\bigr)
\nonumber\\
&\leq& \Biggl(\sum_{j=1}^{n-1}
\frac{\lambda_j(\bbA)-\lambda_{j}(\bbA_k)}{|\lambda
_j(\bbA)-z||\lambda_{j}(\bbA_k)-z|} +\frac{1}{|\lambda_n(\bbA
)-z|} \Biggr)I\bigl(U_n^c
\bigr)
\nonumber\\[-8pt]\\[-8pt]
&\leq& K \Biggl(\sum_{j=1}^{n-1}\bigl(
\lambda_j(\bbA)-\lambda_{j}(\bbA _k)\bigr)+1
\Biggr)I\bigl(U_{n}^c\bigr)\nonumber
\\
&\leq& K\bigl(\lambda_1(\bbA)-\lambda_n(\bbA)+1\bigr)I
\bigl(U_n^c\bigr)\leq K(u_0+3).\nonumber %
\end{eqnarray}
Let $\bbu_j(\bbA_k),j=1,\dots,n-1$ be the eigenvectors corresponding
to the eigenvalues $\lambda_j(\bbA_k),j=1,\dots,n-1$. Then $\sum_{j=1}^{n-1}\frac{\bbu_j(\bbA_k)\bbu_j^T(\bbA_k)}{(\lambda
_{j}(\bbA_k)-z)^2}$ is the spectral decomposition of $\bbD_{k}^{-2}$.
We distinguish two cases:
\begin{enumerate}[(ii)]
\item[(i)] When $z\in V_1=\varsigma_u\cup\{z: |\Im(z)|>(u_0-2)/4\}$, via
(\ref{h2}), we then obtain
\[
|\beta_k|I\bigl(U_n^c\bigr)\leq1/\bigl|\Im(z)\bigr|
\leq\max\bigl\{v_1^{-1},4/(u_0-2)\bigr\}\leq K.
\]
Thus, (\ref{add1}) is true for $z\in V_1$.

\item[(ii)] When $z\in V_2= (\varsigma_l \cup\varsigma_r )\cap
\{z\dvtx  |\Im(z)|<(u_0-2)/4\}$, if $U_n^c$ happens, we have $|\lambda
_{j}(\bbA_k)-\Re(z)|\geq\frac{u_0-2}{2}$ since $\Re(z)=\pm u_0$
for $z\in V_2$. A direct calculation shows
\[
\Re \bigl(\bigl(1+\bbq_k^T\bbD_k^{-2}
\bbq_k\bigr)I\bigl(U_n^c\bigr) \bigr)=1+\sum
_{j=1}^{n-1}\frac{(\lambda_{j}(\bbA_k)-\Re(z))^2-|\Im
(z)|^2}{|\lambda_{j}(\bbA_k)-z|^4}\bigl(
\bbq_k^T\bbu_j(\bbA_k)
\bigr)^2I\bigl(U_n^c\bigr)>1.
\]
Therefore, $|1+\bbq^T\bbD_k^{-2}\bbq|I(U_n^c)$ has a lower bound
which, together with (\ref{s8}), implies (\ref{add1}) is true for
$z\in V_2$.
\end{enumerate}
Since $\varsigma_n=V_1\cup V_2$, we finish the proof of Lemma~\ref{lem5}.
\end{pf}

\begin{lem}\label{lem6}
Let $z\in\varsigma_n$ and $\bar{\mu}_k = \frac{1}{\sqrt {np}}(\mathbf{s}_k^T\mathbf{s}_k-p)-\bbq_k^T\bbD_k^{-1}(z)\bbq
_k+E\frac{1}{np}\tr\bbX
\bbD^{-1}(z)\bbX^T$. The following bounds hold
%
\begin{equation}
\label{add2} E|\bar{\mu}_k|^4\leq K
\frac{\delta^4}{n}+K \biggl(\frac
{1}{n^2}+\frac{n}{p^2}+
\frac{1}{np} \biggr)
\end{equation}
and
%
\begin{equation}
\label{ae2} \bigl|E\bar{\mu}_k^3\bigr|=\mathrm{o}
\bigl(n^{-1}\bigr).
\end{equation}
\end{lem}

\begin{pf} Write
\begin{eqnarray*}
\bar{\mu}_k&=&\frac{1}{\sqrt{np}}\bigl(
\mathbf{s}_k^T\mathbf {s}_k-p\bigr)-\gamma
_{k1}+ \biggl(1+z\sqrt{\frac{p}{n}} \biggr) \biggl(
\frac{1}{n}\tr\bbD ^{-1}(z)-\frac{1}{n}\tr
\bbD_k^{-1}(z) \biggr)
\\
&&{}- \biggl(1+z\sqrt{\frac{p}{n}} \biggr) \biggl(\frac{1}{n}\tr\bbD
^{-1}(z)-E\frac{1}{n}\tr\bbD^{-1}(z) \biggr)+
\frac{1}{\sqrt{np}}
\\
&=&L_1-\gamma_{k1}+L_3+L_4+L_5.
\end{eqnarray*}
When the event $U_n^c$ happens, reviewing the proof of the second
result of (\ref{h10}) and via (\ref{add8}), we also have
\[
E|\gamma_{ks}|^4I\bigl(U_n^c
\bigr)\leq K \biggl(\frac{1}{n^2}+\frac
{n}{p^2}+\frac{1}{np}
\biggr),\qquad m=1,2.
\]
Moreover, by (\ref{s1}) and (\ref{s0})
\[
E|\gamma_{ks}|^4I(U_n)= \mathrm{o}
\bigl(n^{-\ell}\bigr).
\]
It follows that
%
\begin{equation}
\label{add3} E|\gamma_{ks}|^4\leq K \biggl(
\frac{1}{n^2}+\frac{n}{p^2}+\frac
{1}{np} \biggr),\qquad m=1,2.
\end{equation}
Using Lemma~8.10 in Bai and Silverstein \cite{BSbook}, (\ref{s1}),
(\ref{s0}) and (\ref
{s8}) we then have
%
\begin{equation}
\label{y5} E|L_1|^4\leq K\delta^4n^{-1},\qquad
E|L_3|^4\leq Kn^{-4},\qquad  E|L_5|^4
\leq Kn^{-2}p^{-2}.
\end{equation}
As for $L_4$, by Burkholder's inequality, (\ref{u4}) and (\ref{s8}),
we have
%
\begin{eqnarray}
\label{v4} %
E|L_4|^4&\leq&Kn^{-4}E
\Biggl|\sum_{k=1}^{n}(E_k-E_{k-1})
\bigl(\tr\bbD ^{-1}-\tr\bbD_k^{-1}\bigr)
\Biggr|^4
\nonumber\\
&\leq& Kn^{-4}\sum_{k=1}^{n}E \bigl|
\tr\bbD^{-1}(z)-\tr\bbD _k^{-1}(z)
\bigr|^4+Kn^{-1}E \Biggl(\sum_{k=1}^{n}E_k
\bigl|\tr\bbD ^{-1}(z)-\tr\bbD_k^{-1}(z)
\bigr|^2 \Biggr)^2
\nonumber\\
&\leq& Kn^{-4}\sum_{k=1}^{n}E \bigl|
\tr\bbD^{-1}(z)-\tr\bbD _k^{-1}(z)
\bigr|^4I\bigl(U_n^c\bigr)
\\
&&{}+Kn^{-4}E \Biggl(\sum_{k=1}^{n}E_k
\bigl|\tr\bbD^{-1}(z)-\tr\bbD _k^{-1}(z)
\bigr|^2 \Biggr)^2I\bigl(U_n^c\bigr)+
\mathrm{o}\bigl(n^{-\ell}\bigr)
\nonumber\\
&\leq& Kn^{-2}.\nonumber %
\end{eqnarray}
Therefore, the proof of (\ref{add2}) is completed. Also, the analysis
above yields
%
\begin{equation}
\label{e3} \hspace*{-5pt}E|L_1-\gamma_{k1}|^4\leq K
\biggl(\frac{\delta^2}{n}+\frac
{1}{n^2}+\frac{n}{p^2}+\frac{1}{np}
\biggr)\leq K\delta^2n^{-1}, \qquad E|L_3+L_4+L_5|^4
\leq Kn^{-2}.
\end{equation}
It is also easy to verify that, for $z\in\varsigma_n$,
%
\begin{equation}
\label{i15} E\biggl|\frac{1}{\sqrt{np}}\bigl(\mathbf{s}_k^T
\mathbf{s}_k-p\bigr)\biggr|^2\leq Kn^{-1}, \qquad E|
\gamma_{km}|^2\leq Kn^{-1}.
\end{equation}

We proceed to prove (\ref{ae2}). First of all
%
\begin{eqnarray}
\label{i14} %
\bigl|EL_1^3\bigr|=\frac{1}{(np)^{3/2}} \Biggl|E
\Biggl(\sum_{j=1}^{p}\bigl(X_{jk}^2-1
\bigr) \Biggr)^3 \Biggr|=\frac{1}{(np)^{3/2}}\sum
_{j=1}^{p}E\bigl(X_{jk}^2-1
\bigr)^3\leq K\delta^2/n. %
\end{eqnarray}
For $s=1,2$, denoting $\bbM_k^{(s)}=(a_{ij}^{(s)})_{p\times p}$, we
then have
\begin{eqnarray*}
E\gamma_{ks}^3&=&\frac{1}{n^3p^3}E \Biggl(
\sum_{i\neq
j}X_{ik}X_{jk}a_{ij}^{(s)}+
\sum_{i=1}^{n}\bigl(X_{ik}^2-1
\bigr)a_{ii}^{(s)} \Biggr)^3
\\
&=&J_1+J_2+J_3+J_4, %
\end{eqnarray*}
where
\begin{eqnarray*}
J_1&=&\frac{1}{n^3p^3}E \biggl(\sum
_{i\neq j,j\neq t,t\neq
i}X_{ik}^2X_{jk}^2X_{tk}^2a_{ij}^{(s)}a_{jt}^{(s)}a_{ti}^{(s)}
\biggr)+\frac{4}{n^3p^3}E \biggl(\sum_{i\neq
j}X_{ik}^3X_{jk}^3
\bigl(a_{ij}^{(s)}\bigr)^3 \biggr)\triangleq
J_{11}+J_{12},
\\
J_2&=&\frac{1}{n^3p^3}E \Biggl(\sum_{i=1}^{p}
\bigl(X_{ik}^2-1\bigr)^3\bigl(a_{ii}^{(s)}
\bigr)^3 \Biggr),
\\
J_3&=&3\frac{1}{n^3p^3}E \biggl(\sum_{i\neq
j}X_{ik}
\bigl(X_{ik}^2-1\bigr)X_{jk}
\bigl(X_{jk}^2-1\bigr)a_{ij}^{(s)}a_{ii}^{(s)}a_{jj}^{(s)}
\biggr),
\\
J_4&=&3\frac{2}{n^3p^3}E \biggl(\sum_{i\neq
j}X_{ik}^2
\bigl(X_{ik}^2-1\bigr)X_{jk}^2a_{ij}^{(s)}a_{ii}^{(s)}a_{ji}^{(s)}
\biggr). %
\end{eqnarray*}
The inequality (\ref{h7}) can be extended to the range $z\in\varsigma
_n$ by a similar method as that in (\ref{add3}). Therefore,
\begin{eqnarray*}
|J_2|&\leq& K\frac{1}{n^3p^3}p\delta^2\sqrt{np}
\bigl(n^4+n^{2}p\bigr)^{3/4}\leq K
\delta^2n^{-1},
\\
|J_3|&\leq& K\frac{1}{n^3p^3}p^2E\|\mathbf{w}_i
\|^3E\|\mathbf{w}_j\| ^{3}+\mathrm{o}
\bigl(n^{-\ell}\bigr)\leq Kp^{-1}+\mathrm{o}\bigl(n^{-\ell}
\bigr),\qquad  J_4\leq Kp^{-1}+\mathrm{o}\bigl(n^{-\ell}
\bigr),
\end{eqnarray*}
where $\mathbf{w}_j^T$ is the $j$th row of $\bbX_k$.

Consider $J_1$ now. We first note that $J_{12}=\mathrm{O}(p^{-1})$. Split
$J_{12}$ as
\begin{eqnarray*}
J_{12}&=&\frac{1}{n^3p^3}E\tr\bigl(\bbX_k
\bbD_k^{-s}\bbX_k^T
\bigr)^3-\frac
{1}{n^3p^3}E\sum_{i\neq t}a_{ii}^{(s)}a_{it}^{(s)}a_{ti}^{(s)}
\\
&&{}+\frac{1}{n^3p^3}E\sum_{i\neq
j}a_{ij}^{(s)}a_{jj}^{(s)}a_{ji}^{(s)}+
\frac{1}{n^3p^3}E\sum_{i\neq
j}a_{ij}^{(s)}a_{ji}^{(s)}a_{ii}^{(s)}+
\frac{1}{n^3p^3}E\sum_{i=1}^{p}
\bigl(a_{ii}^{(s)}\bigr)^3
\\
&\leq& Kn^{-2}+Kp^{-1}. %
\end{eqnarray*}
Thus, we obtain
%
\begin{equation}
\label{i13} \bigl|E\gamma_{ks}^3\bigr|\leq K\bigl(
\delta^2n^{-1}+p^{-1}\bigr).
\end{equation}
It follows from (\ref{e3}), (\ref{i15}) and (\ref{i13}) that
\begin{eqnarray*}
\bigl|E\bar{\mu}_k^3\bigr|&\leq&\bigl|E(L_1-
\gamma _{k1})^3\bigr|+\bigl|E(L_3+L_4+L_5)^3\bigr|+3\bigl|E(L_1-
\gamma_{k1}) (L_3+L_4+L_5)^2\bigr|
\\
&&{} +3\bigl|E(L_1-\gamma_{k1})^2(L_3+L_4+L_5)\bigr|
\\
&\leq&\bigl|EL_1^3\bigr|+\bigl|E\gamma_{k3}^3\bigr|+3E^{1/2}EL_1^4
\cdot E^{1/2}\gamma _{ks}^2+3E^{1/2}L_1^2
\cdot E^{1/2}\gamma_{ks}^4+Kn^{-3/2}+K
\delta n^{-1}
\\
&=& \mathrm{o}\bigl(n^{-1}\bigr). %
\end{eqnarray*}
The proof of Lemma~\ref{lem6} is completed.
\end{pf}


The following lemma will be used to prove the first result of (\ref
{y11}) and (\ref{m5}) below.

%
\begin{lem}\label{lem7}
For $z\in\varsigma_n$ we have
\[
E\bigl|M_n^{(1)}(z)\bigr|\leq K,
\]
where $M_n^{(1)}(z)=n(m_n(z)-Em_n(z))$.
\end{lem}

\begin{pf} Note that the expression $M_n^{(1)}(z)$ in (\ref{u4})
may not be suitable for $z\in\varsigma_n$, since $\beta_k^{\trr}$ or
even $\beta_k^{\trr}I(U_n^c)$ may be not bounded. For this reason, we
introduce the following notations with the purpose to obtain a similar
expression to (\ref{u4}). Let
\[
\acute{\epsilon}_k=\frac{1}{z+(\afrac{1}{np})E\tr\bbM_k^{(1)}},\qquad  \acute {\mu}_k=
\frac{1}{\sqrt{np}}\bigl(\mathbf{s}_k^T
\mathbf{s}_k-p\bigr)-\gamma _{k1}- \biggl(
\frac{1}{np}\tr\bbM_k^{(1)}-\frac{1}{np}E\tr\bbM
_k^{(1)} \biggr).
\]
Hence
%
\begin{equation}
\label{v7}\beta_k=\acute{\epsilon}_k+
\beta_k\acute {\epsilon}_k\acute{\mu}_k.
\end{equation}
As in (\ref{u4}) and a few lines below it, by (\ref{v7}), we write
\begin{eqnarray*}
M^{(1)}(z)=\sum_{k=1}^{n}(E_k-E_{k-1})
(\acute{\iota}_{k1}+\acute {\iota}_{k2}+\acute{
\iota}_{k3}+\acute{\kappa}_k),
\end{eqnarray*}
where
\begin{eqnarray*}
\acute{\iota}_{k1}(z)&=&- \biggl(1+\frac{1}{np}\tr
\bbM_k^{(2)} \biggr) (\acute{\epsilon}_k)^2
\acute{\mu}_k, \qquad \acute{\iota }_{k2}(z)=-
\gamma_{k1}(\acute{\epsilon}_k)^2\acute{
\mu}_k,
\\
\acute{\iota}_{k3}(z)&=&- \biggl(1+\frac{1}{np}
\bbq_k^T\bbD _k^{-2}(z)
\bbq_k \biggr)\beta_k(\acute{\epsilon}_k)^2
\acute{\mu }_k^2,\qquad  \acute{\kappa}_k=\acute{
\epsilon}_k\gamma_{k2}(z).
\end{eqnarray*}
We next derive the bounds for $\acute{\epsilon}_k$ and the forth
moment of $\acute{\mu}_k$. Since $F_n\xrightarrow{\mathrm{a.s.}} F$ as $n\to
\infty$, we conclude from (\ref{s1}), (\ref{s0}), (\ref{add8}) and
the dominated convergence theorem that, for any fixed positive integer $t$
%
\begin{equation}
\label{j5} E\bigl|m_n(z)-m(z)\bigr|^t\to0.
\end{equation}
By (\ref{a0}), (\ref{s8}) and (\ref{j5}), we then have
\begin{eqnarray*}
E\frac{1}{np}\tr\bbM_k^{(1)}=E \biggl[\biggl(1+z
\sqrt{\frac
{n}{p}}\biggr)m_n(z)-\biggl(1+z\sqrt{
\frac{n}{p}}\biggr)\frac{1}{n}\bigl(\tr\bbD ^{-1}-\tr
\bbD_k^{-1}\bigr)+\frac{n-1}{\sqrt{np}} \biggr]\to m(z).
\end{eqnarray*}
Hence,
%
\begin{equation}
\label{m3} |\acute{\epsilon}_k|=\biggl|\frac{1}{z+m(z)+\mathrm{o}(1)}\biggr|\leq\biggl|
\frac
{2}{z+m(z)}\biggr|\leq2.
\end{equation}
On the other hand, via (\ref{a0}), (\ref{s8}) and (\ref{v4})
\begin{eqnarray*}
E\biggl|\frac{1}{np}\tr\bbM_k^{(1)}-E\frac{1}{np}\tr
\bbM_k^{(1)}\biggr|^4\leq \biggl(1+z\sqrt{
\frac{n}{p}}\biggr)^4n^{-4}E\bigl|\tr\bbD^{-1}-E\tr
\bbD^{-1}\bigr|^4\leq Kn^{-2},
\end{eqnarray*}
and this, together with (\ref{add3}), implies
%
\begin{equation}
\label{m4} E|\acute{\mu}_k|^4\leq K
\frac{\delta^4}{n}+K \biggl(\frac
{1}{n^2}+\frac{n}{p^2}+
\frac{1}{np} \biggr).
\end{equation}
Combining (\ref{m3}), (\ref{m4}), Lemma~\ref{lem5}, (\ref{s1}),
(\ref{s0}), (\ref{add8}) with Burkholder's inequality, we obtain
\[
E\bigl|M_n^{(1)}(z)\bigr|^2\leq K.
\]
The proof of the lemma is completed.
\end{pf}

\section{Convergence of $M_n^{(1)}(z)$}\label{M1}

To prove Proposition~\ref{pro1}, we need to establish (i) the
finite-dimensional convergence and the tightness of $M_n^{(1)}(z)$; (ii)
the convergence of the mean function $EM(z)$.
This section is devoted to the first target. Throughout this section,
we assume that $z\in\mathbb{C}_1$ and $K$ denotes a constant which
may change from line to line and may depend on $v_1$ but is independent
of $n$.

\subsection{Application of central limit theorem for
martingales}\label{martingales}

In order to establish the central limit theorem for the martingale
(\ref{h12}), we have to check the following two conditions:

\begin{con}
[(Lyapunov condition)]\label{cond2} For some $a>2$,
\[
\sum_{k=1}^{n}E_{k-1} \Biggl[
\Biggl|E_{k} \Biggl(\sum_{j=1}^l
a_jE_{k}\bigl(\alpha_k(z_j)
\bigr) \Biggr) \Biggr|^a \Biggr]\xrightarrow{i.p.} 0.
\]
\end{con}

\begin{con}\label{cond1} The covariance
%
\begin{equation}
\label{x11} \Lambda_n(z_1,z_2)\triangleq
\sum_{k=1}^{n}E_{k-1}
\bigl[E_{k}\alpha _k(z_1)\cdot
E_{k}\alpha_k(z_2) \bigr]
\end{equation}
converges in probability to $\Lambda(z_1,z_2)$ whose explicit form
will be given in (\ref{d9}).
\end{con}

Condition \ref{cond2} is satisfied by choosing $a=4$, using Lemma~\ref
{lem1}, and the fact that via (\ref{h2})
\begin{eqnarray*}
\bigl|\alpha_k(z)\bigr|=\biggl\llvert \biggl(1+\frac{1}{np}\tr
\bbM_k^{(2)} \biggr) \bigl(\beta_k^{\trr}
\bigr)^2\eta_k+\beta_k^{\trr}
\gamma_k\biggr\rrvert \leq\frac
{1+v_1^{-2}}{v_1^{2}}|\eta_k|+
\frac{1}{v_1}|\gamma_k|.
\end{eqnarray*}

Consider Condition \ref{cond1} now.
%
Note that
\begin{eqnarray*}
\alpha_k(z)=- \biggl(1+\frac{1}{np}\tr\bbM_k^{(2)}
\biggr) \bigl(\beta _k^{\trr}\bigr)^2
\eta_k-\gamma_k\beta_k^{\trr}=
\frac{\partial}{\partial
z} \bigl(\beta_k^{\trr}\eta_k
\bigr).
\end{eqnarray*}
By the dominated convergence theorem, we have
%
\begin{equation}
\label{h4} \Lambda_n(z_1,z_2)=
\frac{\partial^2}{\partial z_2\,\partial z_1}\sum_{k=1}^{n}E_{k-1}
\bigl[E_{k} \bigl(\beta_k^{\trr}(z_1)
\eta_k(z_1) \bigr)\cdot E_{k} \bigl(
\beta_k^{\trr}(z_2)\eta_k(z_2)
\bigr) \bigr].
\end{equation}
By (\ref{a0}), (\ref{h11}), (\ref{h3}), (\ref{q4}) and the fact
$m_n(z)\xrightarrow{\mathrm{a.s.}}m(z)$, and the dominated convergence theorem
again, for any fixed $t$,
%
\begin{equation}
\label{q1} E\biggl|\frac{1}{np}\tr\bbM_k^{(1)}-
m(z)\biggr|^t\to0,\qquad  E\bigl|\beta_k^{\trr}(z)+m(z)\bigr|^t
\to0, \qquad \mbox{as } n\to\infty.
\end{equation}
Substituting (\ref{q1}) into (\ref{h4}) yields
%
\begin{eqnarray}
\label{h6} %
\Lambda_n(z_1,z_2)&=&
\frac{\partial^2}{\partial z_2\,\partial z_1} \Biggl[m(z_1)m(z_2) \sum
_{k=1}^{n}E_{k-1} \bigl( E_{k}
\eta_k(z_1)\cdot E_{k}\eta_k(z_2)
\bigr)+\mathrm{o}_{\mathrm{i.p.}}(1) \Biggr]\nonumber
\\[-8pt]\\[-8pt]
&=&\frac{\partial^2}{\partial z_2\,\partial z_1} \bigl[m(z_1)m(z_2)\tilde{
\Lambda}_n(z_1,z_2)+\mathrm{o}_{\mathrm{i.p.}}(1)
\bigr].\nonumber %
\end{eqnarray}
By Vitali's theorem (see Titchmarsh \cite{Tit39}, page 168), it
is enough to find
the limit of $\tilde{\Lambda}_n(z_1,z_2)$.
%
To this end, with notation $E_{k}(\bbM_k^{(1)}(z))= (a_{ij}(z)
)_{n\times n}$, write
\[
E_{k}\eta_k(z)=\frac{1}{\sqrt{np}}\sum
_{j=1}^{p}\bigl(X_{jk}^2-1
\bigr)-\frac
{1}{np} \Biggl(\sum_{i\neq j}X_{ik}X_{jk}a_{ij}(z)+
\sum_{i=1}^{p}\bigl(X_{ik}^2-1
\bigr)a_{ii}(z) \Biggr).
\]
By the above formula and independence between $\{X_{ik}\}_{i=1}^{p}$
and $E_{k}(\bbM_k^{(1)})$, a straightforward calculation yields
%
\begin{equation}
\label{o3} E_{k-1} \bigl[E_{k}\eta_k(z_1)
\cdot E_{k}\eta_k(z_2) \bigr] =
\frac{1}{n}E\bigl(X_{11}^2-1\bigr)^2+A_1+A_2+A_3+A_4,
\end{equation}
where
\begin{eqnarray*}
A_1&=&-\frac{1}{np\sqrt{np}}E\bigl(X_{11}^2-1
\bigr)^2\sum_{i=1}^p
a_{ii}(z_1),\qquad  A_2=-\frac{1}{np\sqrt{np}}E
\bigl(X_{11}^2-1\bigr)^2\sum
_{i=1}^pa_{ii}(z_2),
\\
A_3&=&\frac{2}{n^2p^2}\sum_{i\neq j}^pa_{ij}(z_1)a_{ij}(z_2),\qquad
A_4=\frac{1}{n^2p^2}E\bigl(X_{11}^2-1
\bigr)^2\sum_{i=1}^pa_{ii}(z_1)a_{ii}(z_2).
\end{eqnarray*}
%
Note that $a_{ii}(z)$ is precisely $E_ka_{ii}^{(1)}$ in (\ref{h7}).
From (\ref{h7}), we then obtain for $j=1,2,4$
\[
E\Biggl|\sum_{k=1}^{n}A_j\Biggr|
\rightarrow0.
\]
Also, we conclude from (\ref{h7}) that
\begin{eqnarray*}
\sum_{k=1}^{n}A_3=
\frac{2}{n}\sum_{k=1}^{n}
\mathbb{Z}_k-\frac
{2}{n^2p^2}\sum_{k=1}^{n}
\sum_{i=1}^{p}a_{ii}(z_1)a_{ii}(z_2)=
\frac
{2}{n}\sum_{k=1}^{n}
\mathbb{Z}_k+\mathrm{o}_{L_1}(1),
\end{eqnarray*}
where
\[
\mathbb{Z}_k=\frac{1}{np^2}\tr E_{k}
\bbM_k^{(1)}(z_1)\cdot E_{k}\bbM
_k^{(1)}(z_2).
\]
Summarizing the above we see that
%
\begin{equation}
\label{c4} \tilde{\Lambda}_n(z_1,z_2)=
\frac{2}{n}\sum_{k=1}^{n}\mathbb
{Z}_k+\nu_4-1+\mathrm{o}_{L_1}(1).
\end{equation}

\subsection{The asymptotic expression of $\mathbb{Z}_k$}\label{exp}

The goal is to derive an asymptotic expression of $\mathbb{Z}_k$ with
the purpose of obtaining the limit of $\tilde{\Lambda}_n(z_1,z_2)$.

\subsubsection{Decomposition of $\mathbb{Z}_k$}

To evaluate $\mathbb{Z}_k$, we need two different decompositions of
$E_k\bbM_k^{(1)}(z)$. With slight abuse of notation, let $\{\bbe
_i,i=1,\ldots,k-1,k+1,\ldots,n\}$ be the $(n-1)$-dimensional unit
vectors with the $i$th (or $(i-1)$th) element equal to 1
and the remaining equal to 0 according as $i<k$ (or $i>k$). Write $\bbX
_k=\bbX_{ki}+\mathbf{s}_i\bbe_i^T$. Define
%
\begin{eqnarray}
\label{c5} %
\bbD_{ki,r}&=&\bbD_k-
\bbe_i\mathbf{h}_i^T=\frac{1}{\sqrt{np}}\bigl(
\bbX _{ki}^T\bbX_{k}-p\bbI_{(i)}
\bigr)-z\bbI_{n-1},\nonumber
\\
\bbD_{ki}&=&\bbD_k-\bbe_i
\mathbf{h}_i^T-\bbr_i\bbe_i^T=
\frac
{1}{\sqrt {np}}\bigl(\bbX_{ki}^T\bbX_{ki}-p
\bbI_{(i)}\bigr)-z\bbI_{n-1},\nonumber
\\[-8pt]\\[-8pt]
\mathbf{h}_i^T&=&\frac{1}{\sqrt{np}}\mathbf{s}_i^T
\bbX_{ki}+\frac
{1}{\sqrt {np}}\bigl(\mathbf{s}_i^T
\mathbf{s}_i-p\bigr)\bbe_i^T,\qquad
\bbr_i=\frac{1}{\sqrt {np}}\bbX_{ki}^T
\mathbf{s}_i,\nonumber
\\
\zeta_i&=&\frac{1}{1+\vartheta_i},\qquad  \vartheta_i=
\mathbf{h}_i^T\bbD _{ki,r}^{-1}(z)
\bbe_i,\qquad  \bbM_{ki}=\bbX_{ki}
\bbD_{ki}^{-1}(z)\bbX _{ki}^T .\nonumber
\end{eqnarray}
Here $\bbI_{(i)}$ is obtained from $\bbI_{n-1}$ with the $i$th (or
$(i-1)$th) diagonal element replaced by zero if $i<k$ (or $i>k$). With
respect to the above notations we would point out that, for $i<k$ (or
$i>k$), the matrix $\bbX_{ki}$ is obtained from $\bbX_{k}$ with the
entries on the $i$th (or $(i-1)$th) column replaced by zero; $\mathbf
{h}_i^T$ is the $i$th (or $(i-1)$th) row of $\bbA_k$ and $\bbr_i$ is
the $i$th (or $(i-1)$th) column of $\bbA_k$ with the $i$th (or
$(i-1)$th) element replaced by zero. $(\bbX_{ki}^T\bbX_{k}-p\bbI
_{(i)})$ is obtained from $(\bbX_{k}^T\bbX_k-p\bbI_{n-1})$ with the
entries on the $i$th (or $(i-1)$th) row and $i$th (or $(i-1)$th)
column replaced by zero.

The notation defined above may depend on $k$. When we obtain bounds or
limits for them such as $\frac{1}{n}\tr\bbD_{ki}^{-1}$ the results
hold uniformly in $k$.

Observing the structure of the matrices $\bbX_{ki}$ and $\bbD
_{ki}^{-1}$, we have some crucial identities,
%
\begin{equation}
\label{c6} %
\bbX_{ki}\bbe_i=\mathbf{0},\qquad
\bbe_i^T\bbD_{ki,r}^{-1}=\bbe
_i^T\bbD_{ki}^{-1}=-z^{-1}
\bbe_i, %
\end{equation}
where $\mathbf{0}$ is a $p$-dimensional vector with all the elements
equal to 0. By (\ref{c6}) and the frequently used formulas
%
\begin{eqnarray}
\label{c7} %
\bbY^{-1}-\bbW^{-1}&=&-
\bbW^{-1}(\bbY-\bbW)\bbY^{-1},\nonumber \\
 \bigl(\bbY +\mathbf{a}
\mathbf{b}^T\bigr)^{-1}\mathbf{a}&=&\frac{\bbY^{-1}\mathbf
{a}}{1+\mathbf{b}^T\bbY
^{-1}\mathbf{a}},
\\
\mathbf{b}^T\bigl(\bbY+\mathbf{a}\mathbf{b}^T
\bigr)^{-1}&=&\frac{\mathbf
{b}^T\bbY
^{-1}}{1+\mathbf{b}^T\bbY^{-1}\mathbf{a}},\nonumber %
\end{eqnarray}
we have
%
\begin{eqnarray}
\label{i9} %
\bbD_k^{-1}-
\bbD_{ki,r}^{-1}&=&-\zeta_i\bbD_{ki,r}^{-1}
\bbe _i\mathbf{h}_i^T\bbD_{ki,r}^{-1},\nonumber
\\[-8pt]\\[-8pt]
\bbD_{ki,r}^{-1}-\bbD_{ki}^{-1}&=&
\frac{1}{z\sqrt{np}}\bbD _{ki}^{-1}\bbX_{ki}^T
\mathbf{s}_i\bbe_i^T.\nonumber %
\end{eqnarray}

We first claim the following decomposition of $E_k\bbM_k^{(1)}(z)$,
for $i<k$,
%
\begin{eqnarray}
\label{h19} %
E_k\bbM_k^{(1)}(z)&=&E_k
\bbM_{ki}-E_k \biggl(\frac{\zeta_i}{znp}\bbM
_{ki}\mathbf{s}_i\mathbf{s}_i^T
\bbM_{ki} \biggr)+E_k \biggl(\frac{\zeta
_i}{z\sqrt{np}}
\bbM_{ki} \biggr)\mathbf{s}_i\mathbf{s}_i^T
\nonumber\\
&&{} +\mathbf{s}_i\mathbf{s}_i^TE_k
\biggl(\frac{\zeta_i}{z\sqrt {np}}\bbM _{ki} \biggr)-E_k \biggl(
\frac{\zeta_i}{z} \biggr)\mathbf{s}_i\mathbf {s}_i^T
\\
& =& B_1(z)+B_2(z)+B_3(z)+B_4(z)+B_5(z).\nonumber
\end{eqnarray}
Indeed, by the decomposition of $\bbX_k$, write
\begin{eqnarray*}
\bbM_k^{(1)}=\bbX_{ki}\bbD_k^{-1}
\bbX_{ki}^T+\bbX_{ki}\bbD _k^{-1}
\bbe_i\mathbf{s}_i^T+\mathbf{s}_i
\bbe_i^T\bbD_k^{-1}\bbX
_{ki}^T+\mathbf{s}_i\bbe_i^T
\bbD_k^{-1}\bbe_i\mathbf{s}_i^T.
\end{eqnarray*}
Applying (\ref{c5}), (\ref{c6}) and (\ref{i9}), we obtain
\begin{eqnarray*}
\bbX_{ki}\bbD_k^{-1}
\bbX_{ki}^T&=&\bbX_{ki}\bbD_{ki,r}^{-1}
\bbX _{ki}^T-\zeta_i\bbX_{ki}^T
\bbD_{ki,r}^{-1}\bbe_i\mathbf{h}_i^T
\bbD _{ki,r}^{-1}\bbX_{ki}^T
\\
&=&\bbM_{ki}-\frac{\zeta_i}{z\sqrt{np}}\bbM_{ki}
\mathbf{s}_i\cdot \frac{1}{\sqrt{np}}\mathbf{s}_i^T
\bbX_{ki}\bbD_{ki,r}^{-1}\bbX _{ki}^T
\\
&=&\bbM_{ki}-\frac{\zeta_i}{znp}\bbM_{ki}
\mathbf{s}_i\mathbf {s}_i^T
\bbM_{ki}. %
\end{eqnarray*}
Similarly,
\begin{eqnarray*}
\bbX_{ki}\bbD_k^{-1}
\bbe_i\mathbf{s}_i^T&=&\frac{\zeta_i}{z\sqrt {np}}
\bbM_{ki}\mathbf{s}_i\mathbf{s}_i^T,\qquad
\mathbf{s}_i\bbe_i^T\bbD
_k^{-1}\bbX _{ki}^T=
\frac{\zeta_i}{z\sqrt{np}}\mathbf{s}_i\mathbf{s}_i^T
\bbM _{ki},
\\
\mathbf{s}_i\bbe_i^T\bbD_k^{-1}
\bbe_i\mathbf{s}_i^T&=&\zeta _i
\mathbf{s}_i\bbe _i^T\bbD_{ki,r}^{-1}
\bbe_i\mathbf{s}_i^T=-\frac{\zeta_i}{z}
\mathbf {s}_i\mathbf{s}_i^T. %
\end{eqnarray*}
%
Summarizing the above and noting $E_k(\mathbf{s}_i)=\mathbf{s}_i$ for $i<k$
yield (\ref{h19}), as claimed.

On the other hand, write\vspace*{2pt}
\[
\bbD_k=\sum_{i=1(\neq k)}^{n}
\bbe_i\mathbf{h}_i^T-z\bbI_{n-1}.
\]
Multiplying by $\bbD_k^{-1}$ on both sides, we have\vspace*{2pt}
%
\begin{equation}
\label{c8} z\bbD_k^{-1}=-\bbI_{n-1}+\sum
_{i=1(\neq k)}^{n}\bbe_i\mathbf
{h}_i^T\bbD_k^{-1}.
\end{equation}
Therefore, by (\ref{c6}), (\ref{i9}) and the fact that $\bbX_k\bbX
_k^T=\sum_{i\neq k}\mathbf{s}_i\mathbf{s}_i^T$, we have\vspace*{2pt}
%
\begin{eqnarray}
\label{c9} %
zE_{k} \bigl(\bbM_k^{(1)}(z)
\bigr)&=&-E_{k}\bigl(\bbX_k\bbX_k^T
\bigr)+\sum_{i=1(\neq k) }^{n}E_{k-1}\bigl(
\bbX_{k}\bbe_i\mathbf{h}_i^T\bbD
_k^{-1}\bbX _k^T\bigr)\nonumber
\\[2pt]
&=&-E_{k} \Biggl(\sum_{i=1(\neq k)}^{n}
\mathbf{s}_i\mathbf{s}_i^T \Biggr)+\sum
_{i=1(\neq k) }^{n}E_{k}\bigl(
\zeta_i\mathbf{s}_i\mathbf{h}_i^T
\bbD _{ki,r}^{-1}\bigl(\bbX_{ki}^T+
\bbe_i\mathbf{s}_i^T\bigr)\bigr)\nonumber
\\[-7pt]\\[-7pt]
&=&-(n-k)\bbI_{n-1}-\sum_{i<k}
\mathbf{s}_i\mathbf{s}_i^T+\sum
_{i=1(\neq k)
}^{n}E_{k} \biggl(
\frac{\zeta_i}{\sqrt{np}}\mathbf{s}_i\mathbf {s}_i^T
\bbM _{ki} \biggr)\nonumber
\\[2pt]
&&{}+\sum_{i=1(\neq k) }^{n}E_{k} \bigl(
\zeta_i\vartheta_i\mathbf {s}_i
\mathbf{s}_i^T \bigr). \nonumber%
\end{eqnarray}
Consequently, by splitting $E_{k} (\bbM_k^{(1)}(z_2) )$ as in
(\ref{h19}) for $i<k$ and $z_1E_{k} (\bbM_k^{(1)}(z_1) )$ as
in (\ref{c9}), we obtain\vspace*{2pt}
%
\begin{eqnarray}
\label{c11} %
z_1\mathbb{Z}_k&=&
\frac{z_1}{np^2}\tr E_{k}\bbM_k^{(1)}(z_1)
\cdot E_{k}\bbM_k^{(1)}(z_2)
\nonumber\\[-7pt]\\[-7pt]
&=&C_1(z_1,z_2)+C_2(z_1,z_2)+C_3(z_1,z_2)+C_4(z_1,z_2),\nonumber
\end{eqnarray}
where\vspace*{2pt}
\begin{eqnarray*}
C_1(z_1,z_2)&=&-
\frac{1}{np^2}(n-k)\tr E_{k}\bbM_k^{(1)}(z_2),
\\[1pt]
C_2(z_1,z_2)&=&-\frac{1}{np^2}\sum
_{i<k}\mathbf{s}_i^T \Biggl(
\sum_{j=1}^5B_j(z_2)
\Biggr)\mathbf{s}_i=\sum_{j=1}^{5}C_{2j},
\\
C_3(z_1,z_2)&=&\frac{1}{np^2}\sum
_{i<k}E_{k} \Biggl[\frac{\zeta
_i(z_1)}{\sqrt{np}}
\mathbf{s}_i^T\bbM_{ki}(z_1)
\Biggl(\sum_{j=1}^5B_j(z_2)
\Biggr)\mathbf{s}_i \Biggr]
\\
&&{} +\frac{1}{np^2}\sum_{i>k}E_{k}
\biggl[\frac{\zeta_i(z_1)}{\sqrt {np}}\mathbf{s}_i^T
\bbM_{ki}(z_1)E_{k}\bbM_k^{(1)}(z_2)
\mathbf {s}_i \biggr]=\sum_{j=1}^{6}C_{3j},
\\
C_4(z_1,z_2)&=&\frac{1}{np^2}\sum
_{i<k}E_{k} \Biggl[\zeta
_i(z_1)\vartheta_i(z_1)
\mathbf{s}_i^T \Biggl(\sum_{j=1}^5B_j(z_2)
\Biggr)\mathbf{s}_i \Biggr]
\\
&&{} +\frac{1}{np^2}\sum_{i>k}E_{k}
\bigl[\zeta_i(z_1)\vartheta _i(z_1)
\mathbf{s}_i^TE_{k}\bbM_k^{(1)}(z_2)
\mathbf{s}_i \bigr]=\sum_{j=1}^{6}C_{4j},
\end{eqnarray*}
where $C_{2j}$ corresponds to $B_j, j=1,\ldots,5$, for example,
$C_{21}=-\frac{1}{np^2}\sum_{i<k}\mathbf{s}_i^T (B_1(z_2)
)\mathbf{s}_i$, and $C_{3j}$ and $C_{4j}$ are similarly defined. Here both
$C_3(z_1,z_2)$ and $C_4(z_1,z_2)$ are broken up into two parts in terms
of $i>k$ or $i<k$. As will be seen, the terms in (\ref{c11}) tend to 0
in $L_1$, except $C_{25},C_{34},C_{45}$. Next let us demonstrate the details.

\subsubsection{Conclusion of the asymptotic expansion of $\mathbb{Z}_k$}

The purpose is to analyze each term in $C_j(z_1,z_2),j=1,2,3,4$. We
first claim the limits of $\zeta_i,\vartheta_i$ which appear in
$C_j(z_1,z_2)$ for $j=2,3,4$:
%
\begin{equation}
\label{q2} \vartheta_i\xrightarroww{L_4} m(z)/z,\qquad
\zeta_i(z)\xrightarroww{L_4}-zm(z), \qquad \mbox{as } n\to
\infty.
\end{equation}
Indeed, by (\ref{c6}) and (\ref{i9}), we have
%
\begin{equation}
\label{c12} \vartheta_i=\frac{1}{znp}\mathbf{s}_i^T
\bbM_{ki}\mathbf{s}_i-\frac
{1}{z\sqrt{np}}\bigl(
\mathbf{s}_i^T\mathbf{s}_i-p\bigr).
\end{equation}
Replacing $\bbM_k^{(m)}$ in $\gamma_{km}(z)$ by $\bbM_{ki}$, by a
proof similar to that of (\ref{h10}), we have
%
\begin{equation}
\label{q3} E \biggl|\frac{1}{np}\mathbf{s}_i^T
\bbM_{ki}\mathbf{s}_i-\frac
{1}{np}\tr\bbM
_{ki} \biggr|^4\leq K \biggl(\frac{1}{n^2}+
\frac{1}{np} \biggr).
\end{equation}
By (\ref{a0}), we then have $\vartheta_i-\frac{1}{zn}\tr\bbD
_{ki}^{-1}\xrightarroww{L_4}0$. To investigate the distance between
$\tr\bbD_{ki}^{-1}$ and $\tr\bbD_k^{-1}$, let $\dot{\bbA}_{ki}$ be
the matrix constructed from $\bbA_k$ by deleting its $i$th (or
$(i-1)$th) row and $i$th (or $(i-1)$th) column and write $\dot{\bbD
}_{ki}\triangleq\dot{\bbD}_{ki}(z)=\dot{\bbA}_{ki}-z\bbI_{n-2}$
if $i<k$ (or $i>k$). We observe that $\dot{\bbD}_{ki}^{-1}$ can be
obtained from $\bbD_{ki}^{-1}$ by deleting the $i$th (or $(i-1)$th)
row and $i$th (or $(i-1)$th) column if $i<k$ (or $i>k$). Then $\tr
\bbD
_{ki}^{-1}-\tr\dot{\bbD}_{ki}^{-1}=-\frac{1}{z}$. By an identity
similar to (\ref{h11}) and an inequality similar to the bound (\ref
{h3}), we also have $|\tr\bbD_{k}^{-1}-\tr\dot{\bbD}_{ki}|\leq1/v_1$.
Hence $|\tr\bbD_{k}^{-1}-\tr\bbD_{ki}^{-1}|\leq(1/v_1+1/|z|)$. From
(\ref{h11}), we have $|\tr\bbD_{k}^{-1}-\tr\dot{\bbD}|\leq1/v_1$ as
well. As $\frac{1}{n}\tr\bbD^{-1}\xrightarroww{L_t}m(z)$ for any fixed
$t$ by the Helly--Bray theorem and the dominated convergence theorem, we
obtain the first conclusion of (\ref{q2}).

Since the imaginary part of $(z\zeta_i)^{-1}$ is $(\Im(z)+\frac
{1}{np}\Im (\mathbf{s}_i^T\bbM_{ki}\mathbf{s}_i))$ whose
absolute value
is greater than $v_1$, we have $|\zeta_i|\leq|z|/v_1$. Consequently,
via (\ref{q4}), we complete the proof of the second consequence of
(\ref{q2}), as claimed.

Consider $C_1(z_1,z_2)$ first. By (\ref{a0}),
%
\begin{eqnarray}
\label{h18} E\bigl|C_1(z_1,z_2)\bigr|&=&E\biggl
\llvert -\frac{1}{n^2p}(n-k)\tr E_{k}\bbM _k^{(1)}(z_2)
\biggr\rrvert\nonumber\\[-8pt]\\[-8pt]
& \leq&\frac{K}{np^2}n^2p=K\frac{n}{p}\to0.\nonumber
\end{eqnarray}

Before proceeding, we introduce the inequalities for further
simplification in the following. By Lemma~8.10 in Bai and Silverstein
\cite{BSbook} and
(\ref{a0}), for any matrix $\bbB$ independent of $\mathbf{s}_i$,
%
\begin{equation}
\label{h17} E\bigl|\mathbf{s}_i^T\bbM_{ki}\bbB
\mathbf{s}_i\bigr|^2\leq K \bigl(E\bigl|\mathbf
{s}_i^T\bbM _{ki}\bbB\mathbf{s}_i-
\tr\bbM_{ki}\bbB\bigr|^2+KE|\tr\bbM_{ki}\bbB
|^2 \bigr)\leq Kp^2n^2E\|\bbB
\|^2,
\end{equation}
where we also use the fact that, via (\ref{a0}),
\begin{eqnarray*}
|\tr\bbM_{ki}\bbB\overline{\bbB}\overline{
\bbM}_{ki} |&=& \bigl|\tr\bbD_{ki}^{-1/2}
\bbX_{ki}^T\bbB\overline{\bbB}\bbX _{ki}
\overline{\bbD}_{ki}^{-1}\bbX_{ki}^T
\bbX_{ki}\bbD _{ki}^{-1/2} \bigr|
\\
&\leq& n\bigl\|\bbD_{ki}^{-1/2}\bbX_{ki}^T
\bigr\|^2\cdot\|\bbB\|^2\cdot\bigl\|\bbX _{ki}\overline{
\bbD}_{ki}^{-1}\bbX_{ki}^T\bigr\|
\\
&=&n\cdot\|\bbB\|^2\cdot\bigl\|\bbD_{ki}^{-1}
\bbX_{ki}^T\bbX_{ki}\bigr\|^2
\\
&=&n\cdot\|\bbB\|^2\cdot\bigl\|p\bbD_{ki}^{-1}+
\sqrt{np}\bigl(I_{n-1}+z\bbD _{ki}^{-1}\bigr)\bigr\|^2
\\
&\leq& Knp^2\|\bbB\|^2. %
\end{eqnarray*}
For $i>k$, since $E_k\bbM_k$ is independent of $\mathbf{s}_i$, we
similarly have
%
\begin{equation}
\label{q5} E\bigl|\mathbf{s}_i^TE_k
\bbM_{k}\bbB\mathbf{s}_i\bigr|^2\leq
Kn^2p^2.
\end{equation}

Applying Cauchy--Schwarz's inequality, (\ref{h17}) with $\bbB=\bbI
_{n-1}$ and the fact that $|\zeta_i|$ is bounded by $|z|/v_1$ we have
%
\begin{equation}
\label{q6} E|C_{2j}|\leq K\sqrt{\frac{n}{p}}, \qquad j=1,2,3,4.
\end{equation}
Using (\ref{h17}) with $\bbB=E_k\bbM_{ki}(z_2)$ or $\bbB=E_k\bbM
_k$ in (\ref{h17}), we also have
%
\begin{equation}
\label{q7} E|C_{3j}|\leq K\sqrt{\frac{n}{p}},\qquad  j=1,2,3,4.
\end{equation}
By (\ref{q5}), (\ref{q2}) and (\ref{h17}) with $\bbB=\bbI_{n-1}$,
we obtain
%
\begin{equation}
\label{q8} E|C_{4j}|\leq K\frac{n}{p},\qquad  j=1,2,3,4,6.
\end{equation}

Consider $C_{32}$ now. Define $\breve{\zeta}_i $ and $\breve{\bbM
}_{ki}$, the analogues of $\zeta_i(z)$ and $\bbM_{ki}(z)$
respectively, by $ (\mathbf{s}_1,\ldots,\mathbf{s}_k,\breve{\mathbf
{s}}_{k+1},\ldots,\breve{\mathbf{s}}_n )^T$, where $\breve{\mathbf
{s}}_{k+1},\ldots,\breve{\mathbf{s}}_n$ are i.i.d. copies of $\mathbf
{s}_{k+1},\ldots,\mathbf{s}_n$ and independent\vspace*{2pt} of $\mathbf
{s}_1,\ldots,\mathbf{s}_n$.
Then $\breve{\zeta}_i, \breve{\bbM}_{ki}$ have the same properties
as $\zeta_i(z),\bbM_{ki}(z)$, respectively.\vspace*{2pt} Therefore, $|\breve
{\zeta}_i|\leq|z|/v_1$ and $\|\breve{\bbM}_{ki}\|\leq Kp$. Applying
(\ref{h17}) with $\bbB=\breve{\bbM}_{ki}(z_1)$, we have
%
\begin{eqnarray}
\label{c17} %
E|C_{32}|&=&E\biggl\llvert
\frac{1}{np^2}\sum_{i<k}E_{k}E_{k}
\biggl(\frac
{\zeta_i(z_1)}{\sqrt{np}}\mathbf{s}_i^T
\bbM_{ki}(z_1)\frac{\breve
{\zeta}_i(z_2)}{z_2np}\breve{
\bbM}_{ki}(z_2)\mathbf{s}_i\mathbf
{s}_i^T\breve{\bbM}_{ki}(z_2)
\mathbf{s}_i \biggr)\biggr\rrvert\nonumber
\\
&\leq&\frac{K}{n^2p^3\sqrt{np}}\sum_{i<k}E^{\sfrac{1}{2}}
\bigl|\mathbf{s}_i^T\bbM_{ki}(z_1)
\breve{\bbM}_{ki}(z_2)\mathbf {s}_i
\bigr|^2 \cdot E^{\sfrac{1}{2}} \bigl|\mathbf{s}_i^T
\breve{\bbM}_{ki}(z_2)\mathbf {s}_i
\bigr|^2
\\
&\leq& K\sqrt{\frac{n}{p}}.\nonumber %
\end{eqnarray}

Third, consider $C_{25}$. In view of (\ref{q2}), it is straightforward
to check that
%
\begin{equation}
\label{q9} C_{25}=-\frac{k}{n}m(z_2)+
\mathrm{o}_{L_1}(1).
\end{equation}

Further, consider $C_{34}$. By (\ref{q2}) and (\ref{h17}), we have
%
\begin{eqnarray}
\label{ad5} %
C_{34}&=&\frac{1}{np^2}\sum
_{i<k}E_{k} \biggl[\frac{\zeta
_i(z_1)}{\sqrt{np}}
\mathbf{s}_i^T\bbM_{ki}(z_1)B_4(z_2)
\mathbf {s}_i \biggr]\nonumber
\\
&=&\frac{1}{np^2}\sum_{i<k}E_{k}
\biggl[\frac{\zeta_i(z_1)}{\sqrt {np}}\mathbf{s}_i^T
\bbM_{ki}(z_1)E_{k} \biggl(\frac{\zeta
_i(z_2)}{z_2\sqrt{np}}
\bbM_{ki}(z_2) \biggr)\mathbf{s}_i\mathbf
{s}_i^T\mathbf{s}_i \biggr]\nonumber
\\
&=&z_1m(z_1)m(z_2)\frac{1}{n^2p^2}\sum
_{i<k}\mathbf{s}_i^TE_{k}
\bbM _{ki}(z_1)\cdot E_{k}\bbM_{ki}(z_2)
\mathbf{s}_i+\mathrm{o}_{L_1}(1)
\\
&=&z_1m(z_1)m(z_2)\frac{1}{n^2p^2}\sum
_{i<k}\tr \bigl(E_{k}\bbM
_{ki}(z_1)\cdot E_{k}\bbM_{ki}(z_2)
\bigr)+\mathrm{o}_{L_1}(1)\nonumber
\\
&=&z_1m(z_1)m(z_2)\frac{k}{n}
\mathbb{Z}_k+\mathrm{o}_{L_1}(1),\nonumber %
\end{eqnarray}
where the last step uses the fact that via (\ref{h19}), (\ref{h17}),
(\ref{c6}) and a tedious but elementary calculation
\[
\frac{1}{np^2} \bigl|\tr \bigl(E_{k}\bbM_{ki}(z_1)
\cdot E_{k}\bbM _{ki}(z_2) \bigr)-\tr E_{k}
\bigl(\bbX_k\bbD_k^{-1}(z_1)
\bbX_k^T \bigr)\cdot E_{k} \bigl(
\bbX_k\bbD_k^{-1}(z_2)
\bbX_k^T \bigr) \bigr|\leq \frac{K}{n}.
\]

Consider $C_{45}$ finally. By (\ref{q2}), we have
%
\begin{equation}
\label{q10} C_{45}=-m^2(z_1)m(z_2)
\frac{k}{n}+\mathrm{o}_{L_1}(1).
\end{equation}

We conclude from (\ref{c11}), (\ref{h18}), (\ref{q6})--(\ref{q10})
and the fact $m^2(z)+zm(z)+1=0$ that
\begin{eqnarray*}
z_1\mathbb{Z}_k&=&-\frac{k}{n}m(z_2)-
\frac{k}{n}m^2(z_1)m(z_2)+
\frac
{k}{n}z_1m(z_1)m(z_2)
\mathbb{Z}_k+\mathrm{o}_{L_1}(1)
\\
&=&\frac{k}{n}z_1m(z_1)m(z_2)+
\frac{k}{n}z_1m(z_1)m(z_2)\mathbb
{Z}_k+\mathrm{o}_{L_1}(1), %
\end{eqnarray*}
which is equivalent to
%
\begin{equation}
\label{d8} \mathbb{Z}_k=\frac{\sklfrac{k}{n}m(z_1)m(z_2)}{1-\sklfrac
{k}{n}m(z_1)m(z_2)}+
\mathrm{o}_{L_1}(1).
\end{equation}

\subsection{Proof of Condition \texorpdfstring{\protect\ref{cond2}}{5.1}}\label{Condition2}

The equality (\ref{d8}) ensures that
\begin{eqnarray*}
&&\frac{1}{n^2p^2}\sum_{k=1}^{n}\tr E_{k}
\bbM_k^{(1)}(z_1) \cdot E_{k}
\bbM_k^{(1)}(z_2)=\frac{1}{n}\sum
_{k=1}^{n}\mathbb{Z}_k
\\
&&\quad \to\quad \int_{0}^{1}\frac{tm(z_1)m(z_2)}{1-tm(z_1)m(z_2)}\dd t=-1-
\bigl(m(z_1)m(z_2) \bigr)^{-1}\log{
\bigl(1-m(z_1)m(z_2) \bigr)}. %
\end{eqnarray*}
Thus, via (\ref{c4}), we obtain
\begin{eqnarray*}
\label{v6} \tilde{\Lambda}_n(z_1,z_2)
\xrightarrow{\mathrm{i.p.}}\nu_4-3-2 \bigl(m(z_1)m(z_2)
\bigr)^{-1}\log{ \bigl(1-m(z_1)m(z_2) \bigr)}.
\end{eqnarray*}
Consequently, by (\ref{h6})
%
\begin{eqnarray}
\label{d9} %
\Lambda(z_1,z_2)&=&
\frac{\partial^2}{\partial z_1\,\partial z_2} \bigl[(\nu_4-3)m(z_1)m(z_2)-2
\log{ \bigl(1-m(z_1)m(z_2) \bigr)} \bigr]
\nonumber\\[-8pt]\\[-8pt]
&=&m'(z_1)m'(z_2) \bigl[
\nu_4-3+2 \bigl(1-m(z_1)m(z_2)
\bigr)^{-2} \bigr].\nonumber %
\end{eqnarray}

\subsection{Tightness of $M_n^{(1)}(z)$}
\label{tight}

This section is to prove the tightness of $M_n^{(1)}(z)$ for $z\in
\mathbb{C}_1$.
By (\ref{h2}) and Lemma~\ref{lem1},
\[
E\Biggl|\sum_{k=1}^{n}\sum
_{j=1}^l a_jE_{k-1}\bigl(
\alpha_k(z_j)\bigr)\Biggr|^2\leq K\sum
_{k=1}^{n}\sum_{j=1}^l|a_j|^2E\bigl|
\alpha_k(z_j)\bigr|^2\leq K,
\]
which ensures condition (i) of Theorem~12.3 of Billingsley
\cite{bili68}.
Condition (ii) of Theorem~12.3 of Billingsley \cite{bili68}
will be verified by
proving
%
\begin{equation}
\label{f1} %
\frac{E|M_n^{(1)}(z_1)-M_n^{(1)}(z_2)|^2}{|z_1-z_2|^2}\leq K, \qquad z_1,z_2
\in\mathbb{C}_1. %
\end{equation}
We employ the same notations as those in Section~\ref{simplification}. Let
\begin{eqnarray*}
\Upsilon_{k1}&=&\frac{1}{np}\mathbf{s}_k^T
\bbX_k\bbD _k^{-1}(z_1) \bigl(
\bbD_k^{-1}(z_1)+\bbD_k^{-1}(z_2)
\bigr)\bbD_k^{-1}(z_2)\bbX
_k^T\mathbf{s}_k
\\
&&{} -\frac{1}{np}\tr\bbX_k\bbD_k^{-1}(z_1)
\bigl(\bbD _k^{-1}(z_1)+\bbD
_k^{-1}(z_2) \bigr)\bbD_k^{-1}(z_2)
\bbX_k^T,
\\
\Upsilon_{k2}&=&\frac{1}{np} \bigl(\mathbf{s}_k^T
\bbX_k\bbD _k^{-1}(z_2)
\bbD_k^{-1}(z_1)\bbX_k^T
\mathbf{s}_k-\tr\bbX_k\bbD _k^{-1}(z_2)
\bbD_k^{-1}(z_1)\bbX_k^T
\bigr),
\\
d_{k1}(z)&=&\beta_k(z) \biggl(1+\frac{1}{np}
\mathbf{s}_k^T\bbM _k^{(2)}(z)
\mathbf{s}_k \biggr),\\
  d_{k2}(z)&=&1+\frac{1}{np}\tr
\bbM _k^{(2)}(z),
\\
d_{k3}&=&1+\frac{1}{np}\tr\bbX_k
\bbD_k^{-1}(z_2)\bbD _k^{-1}(z_1)
\bbX _k^T,\\
 d_{k4}&=&\frac{1}{np}\tr
\bbX_k\bbD_k^{-1}(z_1) \bigl(\bbD
_k^{-1}(z_1)+\bbD_k^{-1}(z_2)
\bigr)\bbD_k^{-1}(z_2)\bbX_k^T.
\end{eqnarray*}
As in (\ref{u4}), we write
\begin{eqnarray*}
&&M_n^{(1)}(z_1)-M_n^{(1)}(z_2)\\
&&\quad =-
\sum_{k=1}^{n}(E_{k}-E_{k-1})
\bigl(d_{k1}(z_1)-d_{k1}(z_2)\bigr)
\\
&&\quad =-(z_1-z_2)\sum_{k=1}^{n}(E_{k}-E_{k-1})
\bigl[\beta_k(z_1) (\Upsilon _{k1}+d_{k4})-
\beta_k(z_1)d_{k1}(z_2) (
\Upsilon_{k2}+d_{k3}) \bigr]
\\
&&\quad =-(z_1-z_2)\sum_{k=1}^{n}(E_{k}-E_{k-1})\\
&&\hphantom{\quad =-(z_1-z_2)\sum_{k=1}^{n}}{}\times\bigl[(l_1+l_2)+l_3-\beta
_k(z_1)\beta_k(z_2)d_{k2}d_{k3}-
\beta_k(z_1)\beta_k(z_2)d_{k3}
\gamma _k(z_2) \bigr]
\\
&&\quad =-(z_1-z_2)\sum_{k=1}^{n}(E_{k}-E_{k-1})
(l_1+l_2+l_3+l_4+l_5+l_6),
\end{eqnarray*}
where
\begin{eqnarray*}
l_1&=&\beta_k(z_1)\Upsilon_{k1},\qquad
l_2=\beta_k(z_1)\beta _k^{\trr}(z_1)
\eta_k(z_1)d_{k4},
\\
l_3&=&-\beta_k(z_1)\Upsilon_{k2}d_{k1}(z_1),\qquad
l_4=-\beta_k(z_1)\beta _k^{\trr}(z_1)
\eta_k(z_1)\beta_k(z_2)d_{k2}(z_2)d_{k3},
\\
l_5&=&-\beta_k^{\trr}(z_1)
\beta_k(z_2)\beta_k^{\trr}(z_2)
\eta _k(z_2)d_{k2}(z_2)d_{k3},\qquad
l_6=-\beta_k(z_1)\beta_k(z_2)d_{k3}
\gamma_k(z_2).
\end{eqnarray*}
Here the last step uses (\ref{h1}) for $\beta_k(z_1)$ and the facts that
\begin{eqnarray*}
\bbD_k^{-2}(z_1)-\bbD_k^{-2}(z_2)&=&(z_1-z_2)
\bbD_k^{-1}(z_1) \bigl(\bbD_k^{-1}(z_1)+
\bbD_k^{-1}(z_2) \bigr)\bbD_k^{-1}(z_2),
\\
\beta_k(z_1)-\beta_k(z_2)&=&(z_2-z_1)
\beta_k(z_1)\beta_k(z_2)\Upsilon
_{k2}+(z_2-z_1)\beta_k(z_1)
\beta_k(z_2)d_{k3},
\\
(E_{k}-E_{k-1} )\beta_k^{\trr}(z_1)d_{k4}&=&0,\qquad
(E_{k}-E_{k-1} )\beta_k^{\trr}(z_1)
\beta_k^{\trr}(z_2)d_{k2}(z_2)d_{k3}=0.
\end{eqnarray*}
By (\ref{a0}) and Lemma~8.10 in Bai and Silverstein \cite{BSbook},
without any tedious
calculations, one may verify that
\begin{eqnarray*}
\label{ad18} %
\bigl|d_{kj}(z)\bigr|\leq K,\qquad  j=1,2,3,4, \quad \mbox{and}\quad
E|\Upsilon_{kj}|^2\leq Kp^{-1}, \qquad j=1,2.
\end{eqnarray*}
The above inequalities, together with Burkholder's inequality, imply
(\ref{f1}).\vadjust{\goodbreak}

\section{Uniform convergence of $EM_n(z)$}\label{M2}
To finish the proof of Proposition~\ref{pro1}, it remains to derive an
asymptotic expansion of $n(Em_n(z)-m(z))$ for $z\in\mathbb{C}_1$
(defined in Section~\ref{strategy}). In order to unify the proof of
Theorem~\ref{thm2} and Corollary~\ref{cor1}, we derive the asymptotic
expansion of $n(Em_n(z)-m(z))$ under both assumptions $n/p\to0$ and
$n^3/p=\mathrm{O}(1)$ in Proposition~\ref{pro2}. For the purpose of proving
(\ref{y1}), we will prove a stronger result in Proposition~\ref
{pro2}, namely uniform convergence of $n(Em_n(z)-m(z))$ for $z\in
\varsigma_n=\bigcup_{i=l,r,u}\varsigma_i$. For $z$ located in the
wider range $\varsigma_n$, the bounds or limits in Section~\ref{TandS} (e.g., Lemma~\ref{lem1}, (\ref{q1}), (\ref{q2})), cannot be
applied directly. Hence in Section~\ref{lemmas}, we re-establish these
and other useful results. Throughout this section, we assume $z\in
\varsigma_n$ and use the same notations as those in Section~\ref{TandS}.
%

\begin{proposition}\label{pro2}
Suppose that assumption \textup{(c)} is satisfied.
\begin{enumerate}[(ii)]
\item[(i)] Under assumption \textup{(b1)}: $n/p\to0$, we have the
following asymptotic expansion
%
\begin{equation}
\label{x44} n \bigl[Em_n(z)-m(z)-\mathcal{X}_n
\bigl(m(z)\bigr) \bigr]=\mathrm{o}(1),
\end{equation}
uniformly for $z\in\varsigma_n=\bigcup_{i=l,r,u}\varsigma_i$, where
$\mathcal{X}_n(m)$ is defined in \eqref{y4}.
\item[(ii)] Under assumption \textup{(b2)}: $n^3/p=\mathrm{O}(1)$, we have the
following asymptotic expansion
%
\begin{eqnarray}
\label{x4}
&&n \biggl[Em_n(z)-m(z)+\sqrt{\frac{n}{p}}m^4(z)
\bigl(1+m'(z) \bigr) \biggr]\nonumber\\[-8pt]\\[-8pt]
&&\quad =m^3(z)
\bigl(m'(z)+\nu_4-2 \bigr) \bigl(1+m'(z)
\bigr)+\mathrm{o}(1),\nonumber
\end{eqnarray}
uniformly for $z\in\varsigma_n=\bigcup_{i=l,r,u}\varsigma_i$.
\end{enumerate}
\end{proposition}

This, together with \eqref{d9} and the tightness of $M_n^{(1)}(z)$ in
Section~\ref{tight}, implies Proposition~\ref{pro1}. It remains to
prove Proposition~\ref{pro2}. To facilitate statements, let
\begin{eqnarray*}
\omega_n = \frac{1}{n}\sum
_{k=1}^{n}m(z)\beta_k\bar{
\mu}_k,\qquad  \bar {\epsilon}_n=\frac{1}{z+E(\afrac{1}{np})\tr\bbX\bbD^{-1}(z)\bbX^T}.
\end{eqnarray*}
Here, $\omega_n,\bar{\epsilon}_n$ all depend on $z$ and $n$, and
$\bar{\epsilon}_n$ are non-random.

%
\begin{lem}\label{lem4}
Let $z\in\varsigma_n$. We have
\begin{eqnarray*}
nE\omega_n=m^3(z) \bigl(m'(z)+
\nu_4-2 \bigr)+\mathrm{o}(1).
\end{eqnarray*}
\end{lem}

Assuming that Lemma~\ref{lem4} is true for the moment, whose proof is
given in Section~\ref{lem4} below, let us demonstrate how to get
Proposition~\ref{pro2}. By (3.8) in Bai \cite{Bai93},
we obtain
%
\begin{equation}
\label{s2} m_n(z)=\frac{1}{n}\tr\bbD^{-1}(z)=-
\frac{1}{n}\sum_{k=1}^{n}
\beta_k.
\end{equation}
Applying \eqref{q4}, \eqref{s2}, \eqref{a0} and taking the
difference between $\beta_k$ and $\frac{1}{z+m(z)}$, we have
%
\begin{eqnarray}
\label{o5} %
Em_n(z)-m(z) &=& -\frac{1}{n}\sum
_{k=1}^{n}E\beta_k +
\frac
{1}{z+m(z)}\nonumber
\\
&=&E\frac{1}{n}\sum_{k=1}^{n}
\beta_km(z) \biggl[\bar{\mu }_k-\bigl(Em_n(z)-m(z)
\bigr)-\sqrt{\frac{n}{p}}\bigl(1+zEm_n(z)\bigr) \biggr]
\\
&=& E\omega_n + m(z) Em_n(z) \bigl(Em_n(z)-m(z)
\bigr) + \sqrt{\frac
{n}{p}}m(z)Em_n(z) \bigl(1+zEm_n(z)
\bigr).\nonumber %
\end{eqnarray}

\textit{Under assumption} $n/p\to0$: Let $Em_n,m$ respectively,
denote $Em_n(z),m(z)$ to simplify the notations below. By \eqref{q4}
and (\ref{o5}), we have
\begin{eqnarray*}
Em_n-m &=& E\omega_n +
m^2(Em_n-m) + m(Em_n-m)^2 +
\sqrt{\frac
{n}{p}}m(Em_n-m) (1+zm)
\\
&&{} + \sqrt{\frac{n}{p}}m^2(1+zm) + \sqrt{\frac{n}{p}}zm(Em_n-m)^2
+ \sqrt{\frac{n}{p}}zm^2(Em_n-m)
\\
& =& \mathcal{A}(Em_n-m)^2 + (\mathcal{B}+1)
(Em_n-m) + \mathcal{C}_n, %
\end{eqnarray*}
where $\mathcal{A},\mathcal{B}$ are defined in \eqref{y4} and
\[
\mathcal{C}_n = E\omega_n-\sqrt{\frac{n}{p}}m^4.
\]
Rearranging the above equation, we observe that $(Em_n-m)$ satisfies
the equation $\mathcal{A}x^2+\mathcal{B}x+\mathcal{C}_n=0$. Solving
the equation, we obtain
\begin{eqnarray*}
x_{(1)}=\frac{-\mathcal{B}+\sqrt{\mathcal{B}^2-4\mathcal
{A}\mathcal{C}_n}}{2\mathcal{A}},\qquad  x_{(2)}=\frac{-\mathcal{B}-\sqrt {\mathcal{B}^2-4\mathcal{A}\mathcal{C}_n}}{2\mathcal{A}},
\end{eqnarray*}
where $\sqrt{\mathcal{B}^2-4\mathcal{A}\mathcal{C}_n}$ is a complex
number whose imaginary part has the same sign as that of~$\mathcal
{B}$. By the assumption $n/p\to0$ and Lemma~\ref{lem4}, we have
$4\mathcal{A}\mathcal{C}_n\to0$. Then $x_{(1)}=\mathrm{o}(1)$ and $x_{(2)} =
\frac{1-m^2}{m}+\mathrm{o}(1)$. Since $Em_n-m=\mathrm{o}(1)$ by (\ref{b2}), we choose
$Em_n-m = x_{(1)}$. Applying Lemma~\ref{lem4} and the definition of
$\mathcal{X}_n(m)$ in \eqref{y4}, we have
\begin{eqnarray*}
n \bigl[Em_n(z)-m(z)-\mathcal{X}_n
\bigl(m(z)\bigr) \bigr] &=& \frac{-4\mathcal
{A} [nE\omega_n-m^3(z) (m'(z)+\nu_4-2 ) ]}{ 2\mathcal
{A} (\sqrt{\mathcal{B}^2-4\mathcal{A}\mathcal{C}_n}+\sqrt {\mathcal{B}^2-4\mathcal{A}\mathcal{C}} )}
\\
&\to&0. %
\end{eqnarray*}
Hence Proposition~\ref{pro2}(i) is proved.

\textit{Under assumption} $n^3/p=\mathrm{O}(1)$: subtracting $m(z)
Em_n(z)(Em_n(z)-m(z))$ on the both sides of \eqref{o5} and then
dividing $\frac{1}{n}(1-m(z)Em_n(z))$, we have
\begin{eqnarray*}
\label{b2} n\bigl(Em_n(z)-m(z)\bigr) &=&
\frac{nE\omega_n}{1-m(z)Em_n(z)} + \sqrt{\frac
{n^3}{p}}\frac{m(z)Em_n(z)(1+zEm_n(z))}{1-m(z)Em_n(z)}
\\
& =& \frac{m^3(z)}{1-m^2(z)} \bigl(m'(z)+\nu_4-2 \bigr) -
\sqrt{\frac
{n^3}{p}}\frac{m^4(z)}{1-m^2(z)} + \mathrm{o} \biggl(\sqrt{
\frac{n^3}{p}} \biggr),
\nonumber
\end{eqnarray*}
where we use \eqref{j5}, Lemma~\ref{lem4}, \eqref{q4} and the fact
that $m'(z)=\frac{m^2(z)}{1-m^2(z)}$. Proposition~\ref{pro2}(ii) is
proved. Hence, the proof of Proposition~\ref{pro2} is completed. Now
it remains to prove Lemma~\ref{lem4}.

\subsection{Proof of Lemma \texorpdfstring{\protect\ref{lem4}}{6.1}}\label{proofLem4}

From the definitions of $\beta_k, \bar{\epsilon}_n$ and $\bar{\mu
}_k$ (see Lemma~\ref{lem6}), we obtain
%
\begin{equation}
\label{s3} \beta_k=\bar{\epsilon}_n+
\beta_k\bar{\epsilon}_n\bar{\mu}_k.
\end{equation}
By (\ref{s3}), we further write $\beta_k$ as $\beta_k=\bar{\epsilon
}_n+\bar{\epsilon}_n^2\bar{\mu}_k+\bar{\epsilon}_n^3\bar{\mu
}_k^2+\beta_k\bar{\epsilon}_n^3\bar{\mu}_k^3$, which ensures that
%
\begin{eqnarray}
\label{s32} %
nE\omega_n&=&m(z)\bar{
\epsilon}_n\sum_{k=1}^{n}E(
\bar{\mu }_k)+m(z)\bar{\epsilon}_n^2\sum
_{k=1}^{n}E\bigl(\bar{\mu}_k^2
\bigr)\nonumber\\
&&{}+ m(z)\bar{\epsilon}_n^3\sum
_{k=1}^{n}E\bigl(\bar{\mu}_k^3
\bigr)+m(z)\bar {\epsilon}_n^3\sum
_{k=1}^{n}E\bigl(\beta_k\bar{
\mu}_k^4\bigr)
\\
&\triangleq&H_1+H_2+H_3+H_4,\nonumber
\end{eqnarray}
where $H_j,j=1,2,3,4$ are defined in the obvious way. As will be seen,
$H_3$ and $H_4$ are both negligible and the contribution to the limit
of $nE\omega_n$ comes from $H_1$ and $H_2$. Now, we analyze
$H_j,j=1,\dots,4$ one by one.

Consider $H_4$ first. It follows from \eqref{q4} and \eqref{j5} that
%
\begin{equation}
\label{m2} %
\bar{\epsilon}_n=\frac{1}{z+m(z)+\mathrm{o}(1)}=-m(z)+
\mathrm{o}(1). %
\end{equation}
By Lemma~\ref{lem5} and Lemma~\ref{lem6},
\begin{eqnarray*}
E\bigl|\beta_k\bar{\mu}_k^4\bigr|\leq KE\bigl|\bar{
\mu}_k^4\bigr|I\bigl(U_n^c\bigr)+E\bigl|
\beta _k\bar{\mu}_k^4\bigr|I(U_n)
\leq K \biggl(\frac{\delta^4}{n}+\frac
{n}{p^2} \biggr)+\mathrm{o}
\bigl(n^{-\ell}\bigr)\leq K\delta^4n^{-1},
\end{eqnarray*}
which, together with (\ref{m2}), further implies
%
\begin{equation}
\label{s16} H_4=\mathrm{o}(1).
\end{equation}
It follows from Lemma~\ref{lem6} and (\ref{m2}) that
%
\begin{equation}
\label{e5} H_3=\mathrm{o}(1).
\end{equation}

Consider $H_1$ next. We have, via (\ref{a0}) and (\ref{h11}),
%
\begin{eqnarray}
\label{s17} %
H_1&=&m(z)\bar{\epsilon}_n
\sum_{k=1}^{n} \biggl(E\frac{1}{np}\tr
\bbX \bbD^{-1}\bbX^T-E\frac{1}{np}\tr
\bbM_k^{(1)} \biggr)\nonumber
\\
&=& \biggl(1+z\sqrt{\frac{n}{p}} \biggr)m(z)\bar{\epsilon}_n
\frac
{1}{n}\sum_{k=1}^{n}E \bigl(\tr
\bbD^{-1}-\tr\bbD_k^{-1} \bigr)+\sqrt {
\frac{n}{p}}m(z)\bar{\epsilon}_n
\\
&=&- \biggl(1+z\sqrt{\frac{n}{p}} \biggr)m(z)\bar{\epsilon}_n
\frac
{1}{n}\sum_{k=1}^{n}E \biggl[
\beta_k \biggl(1+\frac{1}{np}\mathbf {s}_k^T
\bbM_k^{(2)}\mathbf{s}_k \biggr) \biggr] +
\sqrt{\frac{n}{p}}m(z)\bar{\epsilon}_n. \nonumber%
\end{eqnarray}
%
Applying \eqref{s8}, \eqref{add3} and \eqref{j5}, it is easy to see
\begin{eqnarray*}
1+\frac{1}{np}\mathbf{s}_k^T
\bbM_k^{(2)}\mathbf{s}_k = 1+ \biggl(
\frac
{1}{np}\tr\bbM_k^{(1)} \biggr)' +
\mathrm{o}_{L_4}(1) = 1+m'(z)+\mathrm{o}_{L_4}(1).
\end{eqnarray*}
This, together with \eqref{m2}, Lemma~\ref{lem5} and \eqref{s2},
ensures that
%
\begin{equation}
\label{s30} H_1=-m^2(z) \bigl(1+m'(z)
\bigr)Em_n(z)+\mathrm{o}(1)\to-m^3(z)
\bigl(1+m'(z) \bigr).
\end{equation}

Consider $H_2$ now. By the previous estimation of $E\bar{\mu}_k$
included in $H_1$ we obtain
%
\begin{equation}
\label{add7} E\bar{\mu}_k^2=E(\bar{
\mu}_k-E\bar{\mu}_k)^2+\mathrm{O}
\bigl(n^{-2}\bigr).
\end{equation}
Furthermore a direct calculation yields
%
\begin{equation}
\label{s23} E(\bar{\mu}_k-E\bar{\mu}_k)^2=S_1+S_2,
\end{equation}
where
\begin{eqnarray*}
S_1&=&\frac{1}{n}E \bigl(X_{11}^2-1
\bigr)^2+E\gamma_{k1}^2,\qquad  S_2 =
S_{21} + S_{22},
\\
S_{21}&=&\frac{1}{n^2p^2}E \bigl(\tr\bbM_k^{(1)}-E\tr
\bbM _k^{(1)} \bigr)^2,\\
  S_{22} &=& -
\frac{2}{np\sqrt{np}}E \bigl[\bigl(\mathbf{s}_k^T\mathbf
{s}_k-p\bigr) \bigl(\mathbf{s}_k^T
\bbM_k^{(1)}\mathbf{s}_k-E\tr\bbM
_k^{(1)}\bigr) \bigr].
\end{eqnarray*}
We claim that
%
\begin{equation}
\label{add9} nS_1\to\nu_4-1+2m'(z),\qquad
nS_{21}\to0,\qquad  nS_{22}\to0,\qquad  \mbox{as } n\to \infty.
\end{equation}
Indeed, with notation $\bbM_k^{(1)}=(a_{ij}^{(1)})_{p\times
p},i,j=1,\ldots,p$, as illustrated in \eqref{h7}, we have $\frac
{1}{n^2p^2}\sum_{k=1}^{n}\sum_{i=1}^{p}E|a_{ii}^{(1)}|^2\to0$.
Via this, (\ref{j5}) and (\ref{a0}), a simple calculation yields
\begin{eqnarray*}
\label{x12} %
nE\gamma_{k1}^2&=&
\frac{1}{np^2}E \Biggl(\sum_{i\neq
j}X_{ik}X_{jk}a_{ij}^{(1)}+
\sum_{i=1}^{p}\bigl(X_{ik}^2-1
\bigr)a_{ii}^{(1)} \Biggr)^2
\\
&=&\frac{1}{np^2}E \biggl(\sum_{i\neq j}\sum
_{s\neq
t}X_{ik}X_{jk}X_{sk}X_{tk}a_{ij}^{(1)}a_{st}^{(1)}
\biggr)+\frac
{1}{np^2}\sum_{i=1}^{p}E
\bigl[\bigl(X_{ik}^2-1\bigr)^2
\bigl(a_{ii}^{(1)}\bigr)^2 \bigr]
\\
&=&\frac{2}{np^2}E \biggl(\sum_{i,j}a_{ij}^{(1)}a_{ji}^{(1)}
\biggr)+\mathrm{o}(1)=\frac{2}{np^2}E\tr\bigl(\bbM_k^{(1)}
\bigr)^2+\mathrm{o}(1)
\\
&=&\frac{2}{n}E\tr\bbD_k^{-2}+\mathrm{o}(1)
\rightarrow2m^{\prime}(z). %
\end{eqnarray*}
Since $E|X_{11}^2-1|^2=\nu_4-1$, we have proved the first result of
(\ref{add9}). By Burkholder's inequality, Lemma~\ref{lem7}, (\ref
{a0}), (\ref{s1}) and (\ref{s8})
\begin{equation}
\label{m5} n|S_{21}|=K\biggl(1+z\sqrt{\frac{n}{p}}
\biggr)^2\frac
{1}{n}E\bigl|M^{(1)}(z)\bigr|^2+Kn^{-1}
\leq Kn^{-1}.
\end{equation}
Furthermore,
\begin{eqnarray*}
n|S_{22}|&=&\frac{2}{p\sqrt{np}} \Biggl|E \Biggl(\sum
_{t=1}^{p}\bigl(X_{tk}^2-1
\bigr) \Biggr)\cdot \biggl(\sum_{i,j}X_{ik}X_{jk}a_{ij}^{(1)}
\biggr) \Biggr|
\\
&=&\frac{2}{p\sqrt{np}}\bigl|E\bigl(X_{11}^2-1
\bigr)X_{11}^2\cdot E\tr\bbM _k^{(1)}\bigr|
\leq K\sqrt{\frac{n}{p}}+\mathrm{o}\bigl(n^{-\ell}\bigr)\to0.
\end{eqnarray*}
Therefore, the proof of the second result of (\ref{add9}) is
completed. We then conclude from (\ref{add9}), (\ref{add7}), (\ref
{s23}) and (\ref{m2}) that
%
\begin{equation}
\label{s31} H_2\rightarrow m^3(z)
\bigl(2m^{\prime}(z)+\nu_4-1 \bigr).
\end{equation}


Finally, by (\ref{s32}), (\ref{s16}), (\ref{e5}), (\ref{s30}) and
(\ref{s31}), we obtain
\begin{eqnarray*}
\label{s33} nE\omega_n= m^3(z) \bigl(m{'}(z)+
\nu_4-2 \bigr)+\mathrm{o}(1).
\end{eqnarray*}
Lemma~\ref{lem4} is thus proved. This finishes the proof of
Proposition~\ref{pro1}.

\section{Proof of Proposition \texorpdfstring{\protect\ref{pro3}}{3.2}} \label{proofPro3}
Recall the definition of $U_n$ below Proposition~\ref{pro3} or in
Section~\ref{lemmas}. For $i=l,r$, by Lemma~\ref{lem7}
\begin{eqnarray*}
\label{u5} E\biggl\llvert \int_{\varsigma_i}M_n^{(1)}(z)I
\bigl(U_n^c\bigr)\dd z\biggr\rrvert ^2\leq
\int_{\varsigma_i}E\bigl|M_n^{(1)}(z)\bigr|^2
\dd z\leq K\|\varsigma_i\|\to0, \qquad \mbox{as } n\to\infty,
v_1\to0 .
\end{eqnarray*}
Moreover,
\begin{eqnarray*}
\biggl\llvert \int_{\varsigma_i}EM_n(z)I
\bigl(U_n^c\bigr)\dd z\biggr\rrvert \leq\int
_{\varsigma_i}\bigl|EM_n(z)\bigr|\dd z\to0, \qquad \mbox{as } n\to\infty,
v_1\to0,
\end{eqnarray*}
where the convergence follows from Proposition~\ref{pro2}.

\section{Calculation of the mean and covariance}\label{calculation}

To complete the proof of Theorem~\ref{thm2} and Corollary~\ref{cor1},
it remains to calculate the mean function and covariance function of
$Y(f)$ and $X(f)$. The computation exactly follows Bai and Yao
\cite{BY05} and so
we omit it.

\section*{Acknowledgements}
The authors would like to thank the editor, the associate editor and
the referees' insightful suggestions and comments which significantly
improve the quality and the exposition of the work. Particularly the
current mean correction term
and the calibration of the mean correction term benefit from one of the
referees' suggestions. The second author would like to thank Dr. Lin
Liang-Ching for his help with the real data analysis. This work was
partially supported by the Ministry of Education, Singapore, under
grant \# ARC 14/11.



%

\printhistory

\end{document}